\documentclass[11pt]{article}

\usepackage{amsmath}
\usepackage{amssymb}
\usepackage{latexsym}
\usepackage{amsmath, amsfonts,amssymb, amsthm, euscript,makeidx,color,mathrsfs}

\newtheorem{assumption}{Assumption}
\def\qed{ \ \vrule width.2cm height.2cm depth0cm\smallskip}

\newcommand{\ba}{\begin{array}}
\newcommand{\ea}{\end{array}}
\newcommand{\be}{\begin{equation}}
\newcommand{\ee}{\end{equation}}
\newcommand{\bea}{\begin{eqnarray}}
\newcommand{\eea}{\end{eqnarray}}
\newcommand{\beaa}{\begin{eqnarray*}}
\newcommand{\eeaa}{\end{eqnarray*}}

\def\neg{\negthinspace}

\def\dbE{\mathbb{E}}
\def\dbF{\mathbb{F}}

\def\dbP{\mathbb{P}}
\def\dbR{\mathbb{R}}

\def\a{\alpha}
\def\g{\gamma}
\def\d{\delta}
\def\e{\varepsilon}

\def\l{\lambda}

\def\si{\sigma}

\def\f{\varphi}
\def\th{\theta}

\def\f{\phi}
\def\vf{\varphi}

\def\D{\Delta}
\def\G{\Gamma}

\def\O{\Omega}

\def\G{\Gamma}
\def\D{\Delta}

\def\O{\Omega}
\def\cA{{\cal A}}

\def\cF{{\cal F}}

\def\cH{{\cal H}}

\def\cP{{\cal P}}

\def\hE{\mathbb{E}}
\def\hF{\mathbb{F}}

\def\hL{\mathbb{L}}

\def\hR{\mathbb{R}}

\def\sP{\mathscr{P}}

\def\no{\noindent}

\def\ss{\smallskip}
\def\ms{\medskip}
\def\bs{\bigskip}
\def\q{\quad}
\def\qq{\qquad}

\def\pa{\partial}
\def\cd{\cdot}
\def\cds{\cdots}

\def\td{\nabla}

\def\qed{ \hfill \vrule width.25cm height.25cm depth0cm\smallskip}

\newcommand{\basa}{\begin{assumption}}
\newcommand{\easa}{\end{assumption}}

\newcommand{\ol}{\overline}
\newcommand{\ul}{\underline}

\newcommand{\bas}{\begin{assum}}
\newcommand{\eas}{\end{assum}}

\def\limsup{\mathop{\overline{\rm lim}}}
\def\liminf{\mathop{\underline{\rm lim}}}

\def\pa{\partial}

 \def\cd{\cdot}
\def\cds{\cdots}

\def\ol{\overline}
\def\ul{\underline}

\def\dis{\displaystyle}

\def\1{{\bf 1}}

\def\:{\!:\!}
\def\reff#1{{\rm(\ref{#1})}}
\def \proof{{\noindent \bf Proof\quad}}

at 9pt

\def\prehp(#1,#2){\ensuremath{  #1 \cdot #2 }}

\begin{document}

\newtheorem{thm}{Theorem}[section]
\newtheorem{lem}[thm]{Lemma}
\newtheorem{cor}[thm]{Corollary}
\newtheorem{prop}[thm]{Proposition}
\newtheorem{rem}[thm]{Remark}
\newtheorem{eg}[thm]{Example}
\newtheorem{defn}[thm]{Definition}
\newtheorem{assum}[thm]{Assumption}

\renewcommand {\theequation}{\arabic{section}.\arabic{equation}}
\def\thesection{\arabic{section}}

\title{\bf  Convergence Analysis for Entropy-Regularized Control Problems: A Probabilistic Approach\footnote{An earlier version of this paper was titled "On Convergence Analysis of Policy Iteration Algorithms for Entropy-Regularized Stochastic Control Problems". The authors would like to thank Yufei Zhang for very helpful discussions.}}

\author{
Jin Ma\thanks{\noindent Department of
Mathematics, University of Southern California, Los Angeles, 90089; email: jinma@usc.edu. This author is supported in part by NSF grant  \#DMS-2205972.},
~ ~ Gaozhan Wang\thanks{ \noindent Department of
Mathematics, University of Southern California, Los Angeles, 90089;
email: gaozhanw@usc.edu.  } ~ and ~ Jianfeng Zhang \thanks{ \noindent Department of
Mathematics, University of Southern California, Los Angeles, 90089;
email: jianfenz@usc.edu. This author is supported in part by NSF grant  \#DMS-2205972. }}

\date{\today}
\maketitle

\begin{abstract}
In this paper we investigate the convergence of the {\it Policy Iteration Algorithm} (PIA) for a class of general continuous-time  entropy-regularized stochastic control problems. In particular, instead of employing  sophisticated PDE estimates for the iterative PDEs involved in the algorithm (see, e.g., Huang-Wang-Zhou \cite{Huang2023}), we shall provide a simple proof from scratch for the convergence of the PIA. Our approach builds on probabilistic representation formulae for solutions of PDEs and their derivatives.
Moreover, in the finite horizon model and in the infinite horizon model with large discount factor, the similar arguments lead to a super-exponential rate of convergence without tear. Finally, with some extra efforts we show that our approach can be extended to the diffusion control case in the one dimensional setting, also with a super-exponential rate of convergence. 
 \end{abstract}

\vfill \bs
\no{\bf Keywords.} Reinforcement learning, model uncertainty, entropy-regularization,
policy  iteration algorithm, Feynman-Kac formula, Bismut-Elworthy-Li formula, rate of convergence.

\bs

\no{\it 2020 AMS Mathematics subject classification:} {93E35, 60H30, 35Q93.}

\bs

\eject
\section{Introduction}
\setcounter{equation}{0}

The {\it Policy Iteration Algorithm} (PIA), also known as the {\it Policy Improvement Algorithm},  is a well-known approach in numerical optimal control theory. Its main idea is to construct an iteration scheme for the control actions that traces the maximizers/minimizers of the Hamiltonian, so that the corresponding returns are naturally improving. In the discrete time setting, the algorithm can be traced back to Bellman \cite{Bellman1, Bellman2} and Howard \cite{Howard}, see also Bokanowski-Maroso-Zidani \cite{BMZ}, Puterman-Brumelle \cite{PB}, Santos-Rust \cite{SR}, and the references therein. There have also been many works on continuous time models, see, e.g., Ito-Reisinger-Zhang \cite{IRZ}, Jacka-Mijatovi\'{c} \cite{PIA1}, Kerimkulov-Siska-Szpruch \cite{PIA2, PIA3}, and Puterman \cite{PIA4}. In particular, such an algorithm could converge with a super-linear (or super-exponential) rate, see \cite{BMZ, PB, SR} in the discrete time setting and \cite{IRZ} in the continuous time setting.

Motivated by the model uncertainty, the {\it Reinforcement Learning} (RL) algorithms for the entropy-regularized stochastic control problems have received very strong attention in recent years.  By using the idea of relaxed control, the control problem is regularized (or ``penalized") by Shannon's entropy in order to capture the trade-off between exploitation (to optimize) and exploration (to learn the model). Such entropy-regularized problem in a continuous time model was introduced by Wang-Zariphopoulou-Zhou \cite{WZZ}, see also Guo-Xu-Zariphopoulou \cite{GXZ}, Reisinger-Zhang \cite{RZ}, Sethi-Siska-Zhang \cite{SSZ}, and Tang-Zhang-Zhou \cite{TZZ} in this direction, especially on the relation between the entropy-regularized problem and the original control problem. 

The convergence of the PIA for the entropy-regularized problem, in terms of both the value and the optimal strategy, is clearly a central issue in the theory. In a linear quadratic model, Wang-Zhou \cite{WZ} solved the problem explicitly and  the convergence is immediate. Our paper is mainly inspired by the work Huang-Wang-Zhou \cite{Huang2023}, which established the desired convergence in a general infinite horizon diffusion model with drift controls. Note that the values of the iterative sequence in PIA are by nature increasing (and bounded), so the main issue is to identify its limit with the true value function of the entropy-regularized control problem. By using some sophisticated Sobolev estimates, \cite{Huang2023} established uniform regularity for the iterative value functions, especially uniform bounds for their derivatives, and then derived the convergence by compactness arguments. We also refer to   Bai-Gamage-Ma-Xie \cite{BGMX} and Dong \cite{Dong} for some related works. 

 The original purpose of our paper is to provide a simple proof from scratch for the convergence results that first appeared in \cite{Huang2023}, but without using any heavy PDE machinery. Our proof builds on the Bismut-Elworthy-Li representation formulae \cite{Bismut, EL} for derivatives of functions, see also Ma-Zhang \cite{MaZhang1} and Zhang \cite{Zhangthesis}. These formulae enable us to establish the uniform bounds of the derivatives of the iterative value functions rather easily (see \S\ref{sect-small} Step 1 below),  and then the desired convergence under the $C^2$-norm follows immediately. We emphasize that our approach, albeit simple, may allow us to consider more general situations where uniform PDE estimates may not be available in the existing literature, or the algorithm does not follow any PDE at all, see, e.g., our ongoing work Ma-Wang-Zhang-Zhou \cite{MWZZ}.

 It turns out that, in the case when the discount factor is sufficiently large, our argument can lead to a super-exponential rate of convergence of PIA, which is new in the literature\footnote{An exponential rate is obtained independently by \cite{TWZ}, which we will comment more soon.}. In fact, instead of applying the compactness arguments, in this case the representation formulae can yield the rate of convergence directly. In particular, while the involved derivatives still have uniform bounds, we do not require them for the proof of the convergence here. We would also like to note that, when the discount factor is small, in general it is not reasonable to expect a good rate of convergence, see Remark \ref{rem-diverge} and Example \ref{eg-diverge} below.  

 Our approach also works in the finite horizon case. In that case, the associated PDEs become parabolic, and we obtain the same super-exponential rate of convergence. Technically, we do not need a constraint corresponding to the large discount factor in the infinite horizon case. To the best of our knowledge, this result is new in the literature.
 
 The widely recognized and much more challenging question is the convergence analysis for the case when the diffusion term contains control. In this paper we shall argue that, in an infinite horizon setting when the state process is scalar, with some extra efforts our approach can still lead to the $C^2$-convergence of the value function, which implies the convergence of the optimal strategy. In fact, in a further special case we are able to obtain again a super-exponential rate of convergence. However, we must note that the general higher dimensional case is much more subtle, and our current arguments may be ineffective as they rely heavily on the scalar assumption.

 When finalizing the present paper, we learned the very interesting recent paper Tran-Wang-Zhang \cite{TWZ}. In the base case of infinite horizon model with drift control and sufficiently large discount factor, \cite{TWZ} obtained an exponential rate of convergence. The main ideas are similar, but they used Schauder estimates from PDE literature, while we proved the required estimates from scratch by using the probabilistic representation formulae. It is remarkable that \cite{TWZ} established the convergence for  multi-dimensional  diffusion controlled models when the  control is small in a certain sense and the discount factor is sufficiently large. 
 The key is again  the crucial uniform estimates for the associated iterative fully nonlinear PDEs, in the spirit of Evans-Krylov theorem. We should mention that \cite{PIA3} also studied the diffusion control case, without entry-regularization though. It will be very interesting to combine these ideas to explore general diffusion controlled models.  

At this point we should point out that, although the entropy-regularization is strongly motivated by model uncertainty, the PIA in our paper (also in \cite{Huang2023, TWZ}) relies on the model parameters, and thus is not truly implementable. Building upon the approach in this paper, we shall propose an implementable algorithm and analyze its convergence in a subsequent work \cite{MWZZ}. Another highly relevant issue is the stability of the algorithm, see, e.g., \cite{IRZ, PIA2} and Reisinger-Stockinger-Zhang \cite{RSZ}. We shall address this issue in our future research.

The rest of the paper is organized as follows. In \S2 we formulate the problem and prove the main results in the finite horizon case. In \S3 we prove the main results in the infinite horizon case, first in the case when the discounting factor is large, and then in the general case when the discounting factor is small. In \S4 we investigate the problem with  diffusion control in a scalar setting. 

\ms

{\bf Notations.} To end this section and to facilitate the reader,  we list the following notations that will be used frequently throughout the paper. Let $E$ be a generic Euclidean space, whose usual inner product is denoted by $x\cd y$, for $x, y\in E$, where all $x\in E$ are column vectors. In particular, for $A, B \in \dbR^{d\times d}$, we denote $A : B := \mbox{trace}(A B^\top)$, where $B^\top$ is the transpose of $B$. Moreover, $I_d$ denotes the $d\times d$ identity matrix. For two Euclidean spaces $E_1, E_2$, and $m, k\ge 0$, we denote $C^{m, k}_b([0, T]\times E_1; E_2)$ to be the set of functions $\phi: [0, T]\times E_1 \mapsto E_2$ which is $m$-th order continuously differentiable in $t\in[0,T]$ and $k$-th order continuously differentiable in 
$x\in E_1$,  such that $\phi$ as well as all its derivatives involved above are bounded. Moreover, $C^k_b(E_1; E_2)$ denotes the subspace where $\phi: E_1 \to E_2$ is independent of the temporal variable $t$. 
Furthermore, for  $\phi\in C^2(E_1; E_2)$ and $\psi\in C^{1,2}([0, T]\times E_1; E_2)$, we denote
\bea
\label{norms}
\left.\ba{lll}
\dis \|\phi\|_0 := \sup_{x\in E_1} |\phi(x)|, \q \|\phi\|_2 := \|\phi\|_0 + \|\phi_x\|_0 + \|\phi_{xx}\|_0;\\
 \dis \|\psi\|_0 := \sup_{t\in [0, T]} \|\phi(t, \cd)\|_0, \q \|\psi\|_{1,2} :=  \sup_{t\in [0, T]}\big[\|\psi(t,\cd)\|_2 + \|\psi_t(t,\cd)\|_0\big]. 
\ea\right.
\eea
We use both notations $\pa_x \phi = \phi_x$ for derivatives, whichever is more convenient.

Finally, let $A\subseteq E$ be a domain, and $\sP_0(A)$ the set of probability densities $\pi$ on $A$, namely $\pi: A\to \hR_+$  such that $\int_A \pi(a) da = 1$. For each $\pi\in\sP_0(A)$, we denote its corresponding Shannon's entropy by
\bea
\label{shannon}
\cH(\pi) := - \int_A \pi(a) \ln \pi(a) da.
\eea
Moreover, for $\f \in \hL^1(E_1 \times A; E_2)$ with generic Euclidean spaces $E_1, E_2$, we denote:
\bea
\label{fpi}
\tilde \f(x, \pi) := \int_A \f(x, a) \pi(a) da, \qq x\in E_1.
\eea

\section{The Finite Horizon Case}
\label{sect-finite}
\setcounter{equation}{0}

 We begin our discussion by fixing a finite time horizon $[0, T]$, as in this case we have complete results with simple arguments
 and we  shall consider the infinite horizon case ($T =\infty$) in the next two sections. Let $(\O, \cF, \dbP)$ be a probability space on which is defined a  standard $d$-dimensional Brownian motion $W$,  $\hF := \dbF^W$, and the control set $A$ be a bounded domain with smooth boundary in some Euclidean space, which in particular has finite volume: $0<|A|<\infty$. We note that in this section we consider only drift controls.

Given $(t, x)\in [0, T]\times \dbR^d$, our underlying control problem is as follows: 
\bea
\label{control0}
\left.\ba{lll}
\dis d X^{\a}_s =  b(X^{\a}_s, \a_s)ds + \si(X^{\a}_s) dW_s;\ss\\
\dis u_0(t, x) := \sup_\a \dbE\Big[g(X^{\a}_T) + \int_t^T r(X^{\a}_s, \a_s) ds \big| X^\a_t = x\Big],
\ea\right.
\eea
where $b:  \hR^d \times A \to \hR^d$, $\si: \hR^d\to \hR^{d\times d}$,  $r: \hR^d\times A \to \hR$, $g: \hR^d \to \hR$ are measurable functions, and $\a$ is an appropriate $A$-valued admissible control. Here for notational simplicity we assume $b, \si, r$ are time homogeneous. All the results in this section will remain true when they depend on $t$. It is well known that, under certain technical conditions,  $u_0$ satisfies an HJB equation, and there is a vast literature on numerical methods for $u_0$, provided that the coefficients $b, \si, r, g$ are known. 

Strongly motivated by numerical methods for the above control problem but with model uncertainty, namely when the coefficients $b, \si, r, g$ are unknown,  we consider instead the {\it entropy-regularized  exploratory optimal control problem}. That is, we consider an associated {\it relaxed control problem} regularized by the {\it Shannon's entropy} for the purpose of exploration. More precisely,  let $\cA_T$ denote the set of functions $\pi: [0, T] \times \hR^d \to 
 \sP_0(A)$, and recall \reff{shannon}, \reff{fpi}. Our entropy-regularized  exploratory optimal control problem associated to \reff{control0} takes the form: 
\bea
\label{controlu}
\qq~ \left.\ba{lll}
\dis d X^{\pi}_s = \tilde b(X^{\pi}_s, \pi(s,X^{\pi}_s))ds + \si(X^{\pi}_s) dW_s;\ss\\
\dis J(t,x; \pi)\neg:=\neg \hE\Big[g(X^{\pi}_T)\neg +\neg\int_t^T\neg\neg \big[\tilde r(X^{\pi}_s, \pi(s, X^{\pi}_s)) + \l \cH(\pi(s, X^\pi_s))\big]ds\Big|X^\pi_t = x\Big];\ss\\
\dis u(t, x) := \sup_{\pi\in \cA_T} J(t,x; \pi).

\ea\right.
\eea
Here $\l>0$ is the exogenous ``temperature" parameter 
capturing the trade-off between exploitation and exploration. We remark that $u \to u_0$ when $\l \downarrow 0$, see e.g. \cite{TZZ}. 

In the rest of this section we shall assume:\footnote{When $b, \si, r$ depend on $t$, we require only their continuity in $t$, which is weaker than the requirement for the Schauder's estimate in the standard PDE literature. Our conditions on $x$ here are weaker than those in \cite{Huang2023}, but are stronger than those in  \cite{TWZ}. However, since our main argument is from scratch, it will allow us to consider more general situations where uniform PDE estimates may not be available in the existing literature, or the algorithm does not follow any PDE at all, see, e.g., our ongoing work \cite{MWZZ}. We also refer to Footnote \ref{comparison} below for the diffusion control case. } 
\begin{assum}
\label{assum-finite}
(i) $b, r$ are measurable in $a$; 

(ii)  For each $a\in A$ and for $\f = b(a,\cd),  r(a,\cd), \si, g$, we have $\f\in C^2_b(\dbR^d; \dbR)$ with $\|\f\|_2 \le C_0$, where $C_0>0$ is a constant independent of $a$.

(iii) $\si$ is uniformly non-degenerate: $ [\si \si^\top](x) \ge {1\over C_0} I_{d}$, $x\in\hR^d$.
\end{assum}
Throughout this paper, we shall denote by $C>0$ a generic constant  depending only on $d$, $\l$, $|A|$, and $C_0$, but not on $T$, and it is allowed to vary from line to line. In particular, when the constant does depend  on $T$, we shall denote it as $C_T$. 

Clearly, under Assumption \ref{assum-finite} the SDE for $X^{\pi}$ in \reff{controlu} has a unique weak solution and $u(t, x)$ is well defined and satisfies the following  exploratory HJB equation:
\bea
\label{HJBu}
\left.\ba{c}
\dis u_t+ {1\over 2} \si\si^\top : u_{xx}  + H(x, u_x) =0,\q u(T,x) = g(x),\ss\\
\dis\mbox{where}\q 
 H(x, z) := \sup_{\pi\in \sP_0(A)} \big[~\tilde b(x, \pi) \cd z + \tilde r(x, \pi) + \l \cH(\pi)\big].
\ea\right.
\eea
Moreover, a straightforward calculation along the lines of calculus of variation for the Hamiltonian $H$ shows that the optimal relaxed control $\pi^*$ takes form 
\bea
\label{pi*}
\pi^*(t, x, a):=\G(x, u_x(t, x), a), \q (t, x, a)\in [0, T]\times \hR^d\times A,
\eea
 with $\G$ being the Gibbs function: 
\bea
\label{HGamma}
\qq \G(x, z, a) := \frac{\g(x,z,a)}{\int_A \g(x,z,a') d a'},\q \g(x,z,a):= \exp \Big(\frac{1}{\lambda}[b(x, a) \cdot z+r(x, a)]\Big).
\eea
Consequently, the Hamiltonian $H$ can be written in the following form:
 \bea
 \label{H}
  H(x,z)=\lambda \ln \Big(\int_A \g(x,z,a) d a\Big).
  \eea
  We then have the following simple result. 

\begin{lem}
\label{lem-H}
Let Assumption \ref{assum-finite} hold. 

 (i)  $H$ is twice continuously differentiable in $(x, z)$, with the following estimates:
\bea
\label{Hproperty}
|H_z|\le C, ~0\le H_{zz}\le CI_d;\q \big[|H|+ |H_x|+ |H_{xz}|\big](x,z) \le C[1+|z|].
\eea

(ii)  The PDE \reff{HJBu} has a unique  classical solution $u$  with $\|u\|_{1,2}\le Ce^{CT}$. 
\end{lem}

\ms
\proof  (i) Since $b, r$ are twice continuously differentiable in $x$, by \reff{HGamma}, \reff{H} clearly $H$ is also twice continuously differentiable in $(x, z)$. Moreover, it is easy to check that 
\beaa
&\dis H_z = {\int_A b\g d a\over \int_A \g d a},\q H_{zz} = {\int_A b b^\top \g da \int_A \g da - \int_A b \g da \int_A b^\top \g da \over (\int_A \g da)^2};\\
&\dis H_x ={\int_A [ b_x z +  r_x]\g d a\over \int_A \g d a},\q H_{xz} = {\int_A  b_x \g d a\over \int_A \g d a} -{\int_A [ b_x z +  r_x]\g d a\int_A b^\top \g da\over (\int_A \g d a)^2},
\eeaa
where we suppressed the variables. Then it is straightforward to verify \reff{Hproperty}. 

(ii) Given \reff{Hproperty} and the uniform non-degeneracy of $\si$, this result is standard. We refer to \cite[Theorem 2.4.1]{Zhangthesis} for a probabilistic argument,\footnote{We note that the probabilistic arguments do not involve the additional H\"{o}lder continuity, neither on the coefficients nor on $u$, as typically required in the Schauder's estimate in the PDE literature. } which will be used to prove the convergence in this paper.
\qed

\ms
We now introduce the  {\it Policy Iteration Algorithm}  (PIA) for solving PDE \reff{HJBu}:

\ss
\no{\it Step 0.} Set $\dis u^{0}(t,x) := -C_0 -\big[ C_0 -\l (\ln |A|)^+\big](T-t)$;\footnote{Note that $\sup_{\pi\in \sP_0(A)} \cH(\pi)= (\ln |A|)^+$. We set $u^0$ in this way so that $u^0(t, x) \le J(t, x,\pi)$ for all $\pi$. In particular, while it is not crucial for the remaining analysis, this will imply that $u^0 \le u^1$.}

\ms
\no{\it Step n.}  Define $\pi^{n}(t, x, a):=  \Gamma\left(x,  u_x^{n-1}(t, x), a\right)$ and $u^{n}(t, x):=J(t, x, \pi^n)$, $n\ge 1$.

\ms
\no Then, using (\ref{HGamma}) and (\ref{H}), one shows that 
\beaa
H_z(x,u_x^{n-1})  = \tilde b(x, \pi^n(t, x, \cd))= \int_A b(x, a) \G(x,u_x^{n-1}, a) da,
\eeaa
 and that  $u^{n}$ satisfies the following recursive linear PDE (suppressing variables):
\bea 
\label{PDEtn}
&&u_t^n+ \frac{1}{2} \sigma\sigma^\top :   u_{xx}^{n} +H_z(x,u_x^{n-1})\cd (u_x^{n}\neg-\neg u_x^{n-1} )+ H(x,u_x^{n-1} )=0, 
u^n(T,x) = g(x).
\eea
The following result is standard, in particular part (i) follows from similar arguments as in \cite{Huang2023} which deals with infinite horizon models, we thus omit the proof.

\begin{prop}
\label{prop-un}
Let Assumption \ref{assum-finite} hold. Then

\ss
(i) $u^{n}$ is increasing in $n$ and $u^{n}(t, x) \le C_0 + \big[ C_0 +\l (\ln |A|)^+\big](T-t)$;

\ss
(ii) For $n\ge 1$, $u^{n}\in C^{1,2}_b([0, T]\times \hR^d; \hR)$ is the unique classical solution of \reff{PDEtn}.
\end{prop}

\ms
 This clearly indicates that $u^{n}\uparrow u^*$ for some function $u^*$. Our   purpose is to argue that $u^* = u$ and to obtain the rate of convergence.  Our main result is as follows.

\begin{thm}
\label{thm-finite}
Under Assumption \ref{assum-finite},\footnote{We assume the twice differentiability in  Assumption \ref{assum-finite} in order to get the $C^{1,2}$-convergence of $u^n$ in the theorem. If we content ourselves with the $C^{0,1}$-convergence, from our proofs one can easily see that the second order differentiability is not required. We also note that the $C^{0,1}$-convergence of $u^n$ is sufficient for the convergence of the optimal strategies $\pi^n$. } there exists $0<\eta<1$, depending only on $d$ and $C_0$, such that, denoting $\D u^n:=u^{n} - u$, $\D \pi^n := \pi^n - \pi^*$,
\bea
\label{finite-rate}
\dis \|\D u^{n}\|_{0,1} + \|\D\pi^n \|_0\le C_T\eta^{2^n}, \qq\|\D u^{n}\|_{1,2}  \le C_T 2^{-n},
\eea
\end{thm}

To facilitate the proof for the super-exponential rate of convergence\footnote{Alternatively, this super-exponential rate of convergence of $\D u^n$, also called super-linear rate of convergence, can be expressed as $\dis\lim_{n\to \infty}{ \|\D u^{n}\|_{0,1}\over  \|\D u^{n-1}\|_{0,1}}=0$.} in the first part of \reff{finite-rate}, we first establish a lemma.
\begin{lem}
\label{lem-superexponential}
Assume that a sequence of positive numbers $\{\e_n\}_{n\ge 0}$ satisfy:
\bea
\label{super}
\e_n \le {1\over 2} I(\e_{n-1}) + C_1 \eta_0^{2^n},\q\mbox{where}\q I(\e) := \e\wedge \e^2,
\eea
 for some $0\le\eta_0<1$ and $C_1>0$. Then there exist $0\le\eta_1<1$ and $C_2>0$, depending only on $\eta_0, C_1$, and $\e_0$, such that  $\e_n \le C_2 \eta_1^{2^n}$.
\end{lem}
\proof First, \reff{super} implies $\e_n \le {1\over 2}  \e_{n-1} + C_1 \eta_0^n$. One can show by induction that
\bea
\label{super-en1}
\e_n \le {\e_0\over 2^n} + C_1\sum_{i=1}^n {\eta_0^i \over 2^{n-i}} \le Cn (\eta_0 \vee {1\over 2})^n.
\eea
Next, assuming without loss of generality that $C_1 \ge 4$, then by \reff{super} we have
\bea
\label{super-en2}
\qq \tilde \e_n := {\e_n + C_1 \eta_0^{2^n}\over 2} \le ({\e_{n-1}\over 2})^2 + C_1 \eta_0^{2^n} \le ({\e_{n-1} + C_1 \eta_0^{2^{n-1}}\over 2})^2 = (\tilde \e_{n-1})^2. 
\eea
By \reff{super-en1}  we have $\tilde \e_{n_0} < 1$ for some $n_0\ge 1$. Clearly \reff{super-en2} implies $\tilde \e_n \le (\tilde \e_{n_0})^{2^{n-n_0}}$ for all $n\ge n_0$. Then we obtain the desired estimate with $\eta_1 :=  (\tilde \e_{n_0})^{2^{-n_0}}<1$, and $C_2$ large enough so that $\e_n \le C_2 \eta_1^{2^n}$ holds true for $n\le n_0$ as well.
\qed

\ms
\no{\bf Proof of Theorem \ref{thm-finite}.}  We proceed in five steps. 

{\bf Step 1.} We first provide probabilistic representation formulae for $u, u^n$ and their derivatives, which will be crucial for our estimates.  Fix $(t,x)\in [0, T)\times \dbR^d$ and denote
\bea
\label{Xtx}
X^{t,x}_s = x + \int_t^s \si(X^{t,x}_l) dW_l, \q s\in[t,T]. 
\eea
Let $u^n$ be the solution to  (\ref{PDEtn}), then by standard Feynman-Kac formula we have
\bea
\label{unrep}
&&\dis u(t, x) = \dbE\big[g(X^{t,x}_T) + \int_t^T f(s, X^{t,x}_s) ds\big],\nonumber\\
&&\dis u^{n}(t, x) = \dbE\big[g(X^{t,x}_T) + \int_t^T f^n(s, X^{t,x}_s) ds\big], \\
&&\dis f(t, x):= H^u(t, x),\q f^n(t, x) := H^{n-1}_z \cd [u_x^{n}-u_x^{n-1}](t,x)+ H^{n-1}(t, x),\nonumber\\
&&\dis \mbox{where}~  \phi^u(t, x) := \phi(x, u_x(t,x)),~ \phi^n(t,x) := \phi(x, u^n_x(t,x)),~ \phi = H, H_z, \mbox{etc.}\nonumber
\eea
Next, for any $\phi\in C^2_b(\hR^d; \hR)$,  applying the Bismut-Elworthy-Li formula \cite{Bismut, EL} or the representation formula in \cite{MaZhang1}, and following the arguments in \cite{Zhangthesis},  we have \footnote{\label{footnote-rep}\cite{MaZhang1} provides only the first one in \reff{rep}. In the special case $d=1$, $\si \equiv 1$, we have $X^{x}_t := x + W_t$,  $\td X^{x}_t \equiv 1$, $N^{x}_t = {W_t\over t}$, then it is a direct consequence of the integration by parts formula:
\bea
\label{rep-d=1}
\qq \pa_x \hE[\phi(X^{ x}_t)] = \!\! \int_\dbR \!\!  \phi'(x+ y) {1\over \sqrt{2\pi t}} e^{-{y^2\over 2t}}dy = \!\!  \int_\dbR \!\!  \phi(x+ y) {1\over \sqrt{2\pi t}} e^{-{y^2\over 2t}} {y\over t}dy =  \hE\big[\phi(X^{x}_t) N^{x}_t\big].
\eea
The second formula in \reff{rep} follows directly from the arguments in \cite{MaZhang1, Zhangthesis}  by differentiating the first one with respect to the initial value $x$. We also note that, the $N$ in \cite{MaZhang1, Zhangthesis}  is a row vector, corresponding to the transpose of the $N$ here.}
\bea
\label{rep}
 \pa_x \dbE[\phi(X^{t,x}_s)] &=& \dbE\Big[ (\td X^{t,x}_s)^\top \phi_x(X^{t,x}_s) \Big] = \dbE\Big[\phi(X^{t,x}_s) N^{t,x}_s\Big], \nonumber\\
 \pa_{xx} \dbE[\phi(X^{t,x}_s)] &=& \dbE\Big[(\td X^{t,x}_s)^\top \phi_{xx}(X^{t,x}_s) \td X^{t,x}_s+  \phi_x(X^{t,x}_s) \otimes \td^2 X^{t,x}_s \Big] \\
&=&\dbE\Big[ N^{t,x}_s(\phi_x(X^{t,x}_s))^\top \td X^{t,x}_s + \phi(X^{t,x}_s) \td N^{t,x}_s \Big].\nonumber
\eea
Here the $i$-th column of $\td X$ stands for $\pa_{x_i} X$, $\td^2 X\in \dbR^{d\times d\times d}$ is a tensor with $\td^2_j X \in \dbR^{d\times d}$ standing for $\pa_{x_j} \td X$, and by using the Einstein summation for repeated indices:
\bea
\label{DXtx}
\qq \left.\ba{lll}
\dis \td X^{t,x}_s = I_d +  \int_t^s\si_x^i(X^{t,x}_l) \td X^{t,x}_l dW^i_l, \\
\dis \td^2_j X^{t,x}_s \neg= \neg  \int_t^s \big[  \si_{xx_k}^i(X^{t,x}_l) (\td X^{t,x}_l)^{kj} \td X^{t,x}_s 
+ \si_x^i(X^{t,x}_l) \td^2_j X^{t,x}_s \big]dW^i_l,\ss\\
\dis [\f_x(X_s^{t,x})\otimes \td^2 X_s^{t,x}]_{ij}  :=  \f_{x}(X_s^{t,x}) \cd (\td^2_j X_s^{t,x})^i,
\ea\right.
\eea
where $\si^i$  is the $i$-th column of $\si$,  and $(\td^2_j X)^i$ is the $i$-th column of $\td^2_j X$. Similarly, denoting $\check\si := \si^{-1}$ to be the inverse matrix, we have
 \bea
\label{DNtx}
\left.\ba{lll}
\qq \dis N^{t,x}_s:=\frac{1}{s-t}\int_t^s (\check \sigma( X^{t,x}_l)\td X^{t,x}_l)^\top dW_l;\\
\qq \dis \td_i N^{t,x}_s:=\frac{1}{s-t} \int_t^s \big( (\td X^{t,x}_l)^{ij} \check\sigma_{x_j}(X^{t,x}_l) \td X^{t,x}_l + \check\sigma(X^{t,x}_l) \td^2_i X^{t,x}_l\big)^\top dW_l.
\ea\right.
\eea
Furthermore, one can easily check that
\bea
\label{Ntxest}
\qq \hE\big[|\td X^{t,x}_s|^2+|\td^2 X^{t,x}_s|^2\big] \le Ce^{C(s-t)},~ 
\dis \hE\big[|N^{t,x}_s|^2 + |\td N^{t,x}_s|^2\big]\le {Ce^{C(s-t)}\over s-t}.
\eea
Then we have the representation formulae for the first order derivatives
\bea
\label{utxnrep}
\left.\ba{lll}
\dis u_x(t, x) = \dbE\Big[ (\td X^{t,x}_T)^\top g_x(X^{t,x}_T)+\int_t^T f(s, X^{t,x}_s) N^{t,x}_s ds\Big];\\
\dis u_x^{n}(t, x) = \dbE\Big[ (\td X^{t,x}_T)^\top g_x(X^{t,x}_T)+\int_t^T f^n(s, X^{t,x}_s) N^{t,x}_s ds\Big];
\ea\right.
\eea
and those for the second order derivatives
\bea
\label{utxxnrep}
\qq \left.\ba{lll}
\dis u_{xx}^{n}(t, x) = \dbE\Big[ (\td X^{t,x}_T)^\top g_{xx}(X^{t,x}_T) \td X^{t,x}_T+  g_x(X^{t,x}_T) \otimes\td^2 X^{t,x}_T\\
\dis \qq +\int_t^T \big[N^{t,x}_s  (f_x^n(s, X^{t,x}_s))^\top\td X^{t,x}_s  + f^n(s, X^{t,x}_s) \td N^{t,x}_s  \big]   ds\Big],\\
u_{xx}(t, x) = \dbE\Big[(\td X^{t,x}_T)^\top g_{xx}(X^{t,x}_T)\td X^{t,x}_T +  g_x(X^{t,x}_T)\otimes \td^2 X^{t,x}_T\\
\dis\qq +\int_t^T \big[N^{t,x}_s(f_x(s, X^{t,x}_s))^\top \td X^{t,x}_s  + f(s, X^{t,x}_s) \td N^{t,x}_s  \big]   ds\Big].
\ea\right.
\eea


\ms
{\bf Step 2.} In this step we first assume $T\le \d$, for some $\d>0$ which will be specified later. We shall estimate  
$\e^n_1 := \|\D u_x^n\|_0$, here  the subscript $_1$ stands for the first order derivative $\pa_x$ (the meaning of $\e^n_2$ below is therefore clear). Note that, by \reff{unrep},
\beaa
\D f^n := f^n-f = H^{n-1}_z \cd [u_x^{n}- u_x^{n-1}]+ H^{n-1}-H^u.
\eeaa
Recall \reff{Hproperty} that $|H_z|+|H_{zz}|\le C$ and the function $I$ defined in \reff{super}, one can easily see that
$
\big|H^u -H^{n-1}+ H^{n-1}_z\D u_x^{n-1}\big| \leq CI\big(|\D u_x^{n-1}|\big).
$ Again since $|H_z|\le C$, we have
\bea
\label{fnftest}
|\D f^n(t, x)| \le C\big(\e^n_1 + I(\e^{n-1}_1)\big).
\eea
Now for any $n\ge 1$ and $(t,x)\in [0, T]\times \dbR^d$, by \reff{utxnrep}  we have
\beaa
\dis |\D u^n_x(t, x)|  \!\!\!\! &\le&\!\!\!\! \hE\Big[|\D u^n_x(T, X^{t,x}_T)| |\td X^{t,x}_T| + \int_t^T \big|\D f^n(s,X^{t,x}_s)\big| |N^{t,x}_s| ds\Big]\\
\!\!\!\!&\le&\!\!\!\! C\|\D u^n_x(T, \cd)\|_0 \hE[|\td X^{t,x}_T|] + C\big(\e^n_1 + I(\e^{n-1}_1)\big)\int_t^T  \dbE[|N^{t,x}_s|] ds.
\eeaa
We remark that here obviously $\D u^n_x(T, \cd)\equiv 0$, however, for the sake of arguments later we keep this term.
Then by \reff{Ntxest}, for some constant $C_1>0$  independent of $T$,  
\beaa
 |\D u^n_x(t, x) | \le C_1e^{C_1(T-t)} \Big[\|\D u^n_x(T, \cd)\|_0 +\big(\e^n_1 + I(\e^{n-1}_1)\big)\sqrt{T-t}\Big].
 \eeaa
 Since $(t,x)$ is arbitrary and $T-t\le T\le \d$, we obtain that
\bea
\label{en1}
\e^n_1 \le C_1e^{C_1\d} \Big[\|\D u^n_x(T, \cd)\|_0 +\big(\e^n_1 + I(\e^{n-1}_1)\big)\sqrt{\d}\Big].
\eea
We now set $\d>0$ small such that
\bea
\label{delta}
C_1e^{C_1\d } \sqrt{\d} \le {1\over 3}.
\eea
Then, for $T\le \d$, (\ref{en1}) implies $ \e^n_1 \le {1\over 3}\big(\e^n_1 + I(\e^{n-1}_1)\big) + C \|\D u^n_x(T, \cd)\|_0$,
and thus
\bea
\label{en2}
 \e^n_1 \le {1\over 2}  I(\e^{n-1}_1) + C \|\D u^n_x(T, \cd)\|_0.
\eea
Note that $\D u^n_x(T, \cd) =0$. Then $\{\e^n_1\}_{n\ge 0}$ satisfies \reff{super} with $\eta_0=0$, and thus it follows from Lemma \ref{lem-superexponential} that $\e^n_1 \le C\eta^{2^n}$ for some $0<\eta<1$.

{\bf Step 3.} We next estimate $\e^n_1$ for general $T$. Let $0=t_0<\cds<t_m=T$ be such that $t_i - t_{i-1}\le \d$, where $\d$ satisfies \reff{delta} and is independent of $T$.  For each $i$, denote $\e^n_{1,i} :=  \sup_{t \in [t_{i-1}, t_i]} \|\D u^n_x(t,\cd)\|_0$. Apply the arguments in Step 2 on $[t_{i-1}, t_i]$,  by \reff{en2} we have
\bea
\label{en2-i}
 \e^n_{1,i} \le {1\over 2}  I(\e^{n-1}_{1,i}) + C \|\D u^n_x(t_i, \cd)\|_0.
\eea
Note that $\|\D u^n_x(t_m,\cd)\|_0=0 = C_m\eta_m^{2^n}$ with $\eta_m := 0$. Assume 
by induction that $\|\D u^n_x(t_i,\cd)\|_0\le C_i\eta_i^{2^n}$ for some $0<\eta_i <1$. Then, by Lemma \ref{lem-superexponential} we have $\e^n_{1,i}\le C_{i-1}\eta_{i-1}^{2^n}$ for some $0<\eta_{i-1}<1$. In particular, $\|\D u^n_x(t_{i-1},\cd)\|_0\le C_{i-1}\eta_{i-1}^{2^n}$. By a backward induction on $i=m, \cds, 1$, we obtain immediately that
 \bea
 \label{uxnrate}
 \|\D u^n_x\|_0 \le C_T\eta^{2^n},\q\mbox{for some}~ 0<\eta<1.
 \eea

Moreover, by \reff{unrep}, \reff{fnftest}, and \reff{uxnrate}, we have
\beaa
|\D u^n(t,x)| \le \dbE\Big[\int_t^T \big|\D f^n(s, X^{t,x}_s)\big| ds\Big] \le (T-t)\|\D f^n\|_0\le  C_T\eta^{2^n}.
\eeaa
This implies $\|\D u^n\|_0  \le C_T\eta^{2^n}$, and together with \reff{uxnrate}, we have
\bea
\label{unrate}
\|\D u^n\|_{0,1}  \le C_T\eta^{2^n}.
\eea

Furthermore, by \reff{uxnrate} and Lemma \ref{lem-H} (ii) it is clear that $\|u^n_x\|_0, \|u_x\|_0\le C_T$. Then, by \reff{HGamma} one can easily check that, at arbitrary $(t,x,a)$, 
\bea
\label{pirate}
\qq |\D\pi^n| =\big| \Gamma\left(x,  u_x^{n-1}, a\right)-\Gamma\left(x,  u_x, a\right)\big|\le C_T |\D u_x^{n-1}(t, x)| \le C_T \eta^{2^n}.
\eea
This, together with \reff{unrate}, proves the first inequality in \reff{finite-rate}.

 {\bf Step 4.} We now estimate  $\e^n_2 := \|\D u_{xx}^n\|_0$. Again we first assume $T\le \d'$ for some $\d'>0$ to be specified later.  By \reff{utxxnrep} we have
\bea
\label{paxxun-u}
&&\dis |\D u_{xx}^n(t, x)|
 \le \hE\Big[\|\D u_{xx}^n(T,\cd)\|_0 |\td X^{t,x}_T|^2 +  \|\D u_x^n(T,\cd)\|_0  |\td^2 X^{t,x}_T| \nonumber\\
&&\dis\qq  + \int_t^T \big[ N^{t,x}_s(\D f_x^n(s, X^{t,x}_s))^\top \td X^{t,x}_s   + \D f^n(s, X^{t,x}_s) \td N^{t,x}_s  \big]   ds\Big].
\eea
Here again we keep the terminal difference term for our argument. Note that
\beaa
\D f_x^n&=& [H^{n-1}_{xz} + H^{n-1}_{zz} u_{xx}^{n-1}][u_x^n-u_x^{n-1} ]+  H^{n-1}_z\big[u_{xx}^n-u_{xx}^{n-1}\big] \\
 &&+ \big[H^{n-1}_x - H^u_x\big]+\big[H^{n-1}_z u_{xx}^{n-1} - H^u_z u_{xx}\big].
\eeaa
Denote $\hat\e^n_i := \e^n_i+\e^{n-1}_i$, $i=1,2$. Then it follows from  Lemma \ref{lem-H} that 
\bea
\label{paxfn-f}
 \q \big|\D f_x^n\big| \le C\big[(1+\e^{n-1}_2)\hat\e^n_1 +  \hat\e^n_2 +  \e^{n-1}_1 +  \e^{n-1}_2\big]\le   C(1+\hat\e^n_1)\hat\e^n_2 + C\hat\e^n_1.
\eea
 Recall \reff{fnftest} and \reff{uxnrate}, by the arbitrariness of $(t,x)$ we see that \reff{paxxun-u} leads to that
\bea
\label{en2est}
\qq~ ~\e^n_2 \le Ce^{C\d'}\big[\|\D u_{xx}^n(T,\cd)\|_0  +  \| \D u_x^n(T,\cd)\|_0\big] +C_1e^{C_1\d'}\sqrt{\d'}(1+\hat\e^n_1)\hat\e^n_2 + C\hat\e^n_1.
\eea

 {\bf Step 5.} Finally we consider arbitrary $T$ and complete the estimate for $\e^n_2$. Let $0=t_0'<\cds<t_{m'}'=T$ be another partition with $t_i' - t_{i-1}'\le \d'$,  and denote $\e^n_{2,i}:= \sup_{t \in [t_{i-1}', t_i']} \|\D u_{xx}^n(t,\cd)\|_0$. Assume, for the $C_1$ in \reff{en2est} and $C_T, \eta$ in \reff{uxnrate},
  \beaa
 C_1 e^{C_1 \d'}\sqrt{\d'} [1+ C_T\eta^{2^n} + C_T\eta^{2^{n-1}}] \le 1\slash 3.
 \eeaa
 We remark that here we allow $\d'$ to depend on $T$. For each $i$, apply the arguments in Step 4 on $[t_{i-1}', t_i']$, then  by \reff{en2est} and \reff{uxnrate} we have
 \beaa
 \e^n_{2,i} \le  C\|\D u_{xx}^n(t_i',\cd)\|_0   +{1\over 3}\big[\e^n_{2,i} + \e^{n-1}_{2,i}\big] + C_T\eta^{2^{n-1}}.
\eeaa
This implies that
 \beaa
 \e^n_{2,i} \le  {1\over 2}  \e^{n-1}_{2,i} + C\|\D u_{xx}^n(t_i',\cd)\|_0  +C_T\eta^{2^{n-1}}.
\eeaa
We emphasize that the right side here does not involve $(\e^{n-1}_{2,i})^2$. Then by standard arguments one can only obtain an exponential rate of convergence: 
 \beaa
\sup_{t \in [t_{i-1}', t_i']} \|\D u_{xx}^n(t,\cd)\|_0 \le  C\|\D u_{xx}^n(t_i',\cd)\|_0  + C 2^{-n}+ C_T\eta^{2^{n-1}}.
 \eeaa
 Note that $\|\D u_{xx}^n(t_{m'}',\cd)\|_0=0$. Then by backward induction on $i=m', \cds, 1$, we obtain immediately that
 \bea
 \label{uxnrate3}
 \|\D u_{xx}^n\|_0 \le C_T2^{-n}.
 \eea
 
 Finally, by the PDEs \reff{HJBu} and \reff{PDEtn}, we obtain easily from \reff{unrate} and \reff{uxnrate3} the desired estimate for $\|\D u_t^n\|_0$, and thus prove \reff{finite-rate}.
  \qed

\begin{rem}[Dependence of the estimates on $\l$]
\label{rem-lambda} 
In Lemma \ref{lem-H} one can easily see that the bounds of $H_x, H_z, H_{xz}, H_{zz}$ and $u_x, u_{xx}$ are independent of $\l$, but those of $H, u, u_t$ may depend on it. Consequently, by examining the estimates in the above theorem carefully, one can see that the constant $C_T$ for the estimates of $\|\D u\|_{0,1}$ and $\|\D u\|_{1,2}$ in \reff{finite-rate} are independent of $\l$. However, by \reff{HGamma} and \reff{pirate}, the constant $C_T$ for the estimate of $\|\D\pi^n\|$ in \reff{finite-rate} grows exponentially in ${1\over \l}$ and thus explodes when $\l\downarrow 0$. So there is a tradeoff between the numerical efficiency for exploitation (which requires $\l$ to be large) and limiting the error due to the exploration (which requires $\l$ to be small). 
\end{rem}

\section{The Infinite Horizon Case}
\label{sect-infinite}
\setcounter{equation}{0}

In this section we  consider the infinite horizon case: $T=\infty$. We note that in this case the PDEs involved will become purely elliptic, but defined on the whole space. Our argument will be slightly different, but along the same line. We shall first still consider only the drift control case, and will extend our result to the diffusion control case in one dimensional setting in the next section. 

Consider the following entropy-regularized  exploratory optimal control problem: 
\bea
\label{control}
\left.\ba{lll}
\dis d X^\pi_t = \tilde b(X^\pi_t, \pi(t,X^\pi_t))dt +  \si(X^\pi_t) dW_t;\ss\\
\dis J(x, \pi):= \dbE\Big[\int_0^\infty e^{-\rho t} \big[\tilde r(X^\pi_t, \pi(t, X^\pi_t)) + \l \cH(\pi(t, X^\pi_t))\big]dt\Big|X^\pi_0=x\Big], \ss\\
\dis v(x) := \sup_{\pi\in \cA} J(x, \pi),
\ea\right.
\eea
where $\rho>0$ is the discount factor, and $\cA$ denote the set of $\pi: [0, \infty) \times \dbR^d \to 
 \sP_0(A)$. We shall assume that the coefficients $b, \si, r$ satisfy Assumption \ref{assum-finite} with the constant $C_0$.
Throughout this section, the generic constant $C>0$ does not depend on $\rho$. In particular, when the constant does depend  on $\rho$, we shall denote it as $C_\rho$. 

Clearly, for the same Gibbs form $\G$ and function $H$ as in (\ref{HGamma}), (\ref{H}), $v$ is well defined and satisfies the following HJB equation:
\bea
\label{HJBv}
\dis \rho v(x) = {1\over 2} [\si\si^\top](x) : v_{xx} (x)  + H(x,  v_x),\q x\in\hR^d.
\eea
Similarly to Lemma \ref{lem-H}, it is standard to show that, under Assumption \ref{assum-finite} the PDE \reff{HJBv} has a unique  classical solution  $v\in C^2_b(\hR^d;\hR)$, with $\|v\|_2\le C_\rho$. Moreover, the optimal control of \reff{control} is $\pi^*(x, a) = \G(x, v_x(x),a)$.

\ms
Consider now the  {\it Policy Iteration Algorithm}  (PIA) for solving \reff{HJBv} recursively:

\ms
{\it Step 0.} Set $\dis v^{0} := -{1\over \rho}\big[ C_0 -\l (\ln |A|)^+\big]$;

\ms
{\it Step n.}  For $n\ge 1$, define $\pi^{n}(x, a):=  \Gamma\left(x,  v_x^{n-1}(x), a\right)$ and $v^{n}(x):=J(x, \pi^n)$.

\ms
\no Similarly to the analysis in the last section, we can easily check that each $v^{n}$ satisfies the following recursive linear PDE:
\bea 
\label{PDEn}
\qq \rho v^{n}= \frac{1}{2} \sigma\sigma^\top :   v_{xx}^{n} + H_z(x, v_x^{n-1}) \cd ( v_x^{n}\neg-\neg  v_x^{n-1} ) + H(x, v_x^{n-1} ).
\eea
Then we have the following analogue of Proposition \ref{prop-un}, whose proof is omitted.

\ss
\begin{prop}
\label{prop-vn}
Let Assumption \ref{assum-finite} hold. Then

\ss
(i) $v^{n}$ is increasing in $n$ and $v^{n} \le {1\over \rho}\big[ C_0 +\l (\ln |A|)^+\big]$;

\ss
(ii) For each $n\ge 1$, $v^{n}\in C^{2}_b(\dbR^d; \dbR)$ is a classical solution of \reff{PDEn}.
\end{prop}
Our main result of this section is the following analogue of Theorem \ref{thm-finite}, 
\begin{thm}
\label{thm-main1}
Let Assumption \ref{assum-finite} hold and denote $\D v^n:= v^n-v$.

\ss
(i) There exist $\rho_0>0$ and $0<\eta<1$, depending only on  $d$, $C_0$, such that,
\bea
\label{conv2}
 \| \D v^{n}\|_{1}  + \|\D\pi^n\|_0\le C\eta^{2^n}, \q \|\D v^{n}\|_{2}  \le  C2^{-n},\q \mbox{whenever $\rho\ge \rho_0$}.
\eea

(ii) In the general case, $v^n\to v$ in $C^2$ and $\pi^n\to \pi^*$ in $C^0$, uniformly on compacts. That is, for any compact set $K\Subset \hR^d$, it holds that
\bea
\label{conv1}
\lim_{n\to\infty} \|\D v^n\|_{2,K}=0,\q \lim_{n\to\infty} \sup_{a\in A}\|\D\pi^n(\cd, a)\|_{0,K}=0
\eea
where $\|\vf\|_{0, K}:=\sup_{x\in K}|\vf(x)|$, and $\|\vf\|_{2,K}:=\|\vf\|_{0,K} +\|\vf_x\|_{0,K} + \|\vf_{xx}\|_{0,K}$.
\end{thm}

Since the arguments for the proof of Theorem \ref{thm-main1} will depend crucially on the ``size" of the discounting factor $\rho$, we shall carry it out separately in the two subsections below, for ``large" and ``small" $\rho$, respectively.

\ss
\subsection{Proof of Theorem \ref{thm-main1} (i)}
\label{sect-large}

In this subsection we prove \reff{conv2}, assuming that $\rho$ is sufficiently large. We emphasize again that in this subsection the generic constant $C$ does not depend on $\rho$.  We proceed in three steps. 

{\bf Step 1.}  We begin by recalling the  probabilistic representation formulae for $v, v^n$ and their derivatives, which are crucial for our arguments. Denote
\beaa
(X^x, \td X^x, N^x, \td^2_j X^x, \td_i N^x) := (X^{0,x}, \td X^{0,x}, N^{0,x}, \td^2_j X^{0,x}, \td_i N^{0,x}).
\eeaa
Then, similarly to \reff{unrep}, we derive from \reff{HJBv} and \reff{PDEn} that
\bea
\label{vnrep}
&&\dis v(x) = \hE\Big[\int_0^\infty e^{-\rho t} f(X^x_t) dt\Big],\q v^{n}(x) = \hE\Big[\int_0^\infty e^{-\rho t} f^n(X^x_t) dt\Big],\nonumber\\
&&\dis f:= H^v,\q f^n := H^{n-1}_z \cd [ v_x^{n}-  v_x^{n-1} ]+ H^{n-1},\\
&&\dis \mbox{where}\q \phi^v(x) := \phi(x, v_x(x)),~ \phi^n(x) := \phi(x, v^n_x(x)),~ \phi = H, H_z, \mbox{etc.}\nonumber
\eea
Next, applying \reff{rep} on above, we obtain 
\bea
\label{vxnrep}
\qq \left.\ba{lll}
\dis  v_x(x) = \dbE\Big[\int_0^\infty e^{-\rho t} f(X^x_t) N^x_t dt\Big],~   v_x^{n}(x) = \dbE\Big[\int_0^\infty e^{-\rho t} f^n(X^x_t) N^x_t dt\Big];\ss\\
\dis  v_{xx}(x) = \dbE\Big[\int_0^\infty e^{-\rho t}\big[ N^x_t(f_x(X^x_t))^\top \td X^x_t     + f (X^x_t) \td N^x_t\big]   dt\Big];\ss\\
 \dis  v_{xx}^{n}(x) = \dbE\Big[\int_0^\infty e^{-\rho t}\big[N^x_t(f^n_x(X^x_t))^\top \td X^x_t    + f^n(X^x_t) \td N^x_t\big]   dt\Big].
\ea\right.
\eea

{\bf Step 2.} We first estimate $\e^n_1 := \| \D v_x^n\|_0 $. Similarly to \reff{fnftest} we have
\bea
\label{fnfest}
|\D f^n(x)| \le C(\e^n_1 + I(\e^{n-1}_1)),
\eea
for the $I$ in \reff{super}. Then, for any $n\ge 1$ and $x\in \dbR^d$, by \reff{vxnrep}  and \reff{Ntxest} we have
\beaa
&&\!\!\!\!\!\!\! |\D v_x^n(x)| \le  \dbE\Big[\int_0^\infty  e^{-\rho t} \big|\D f^n(X^x_t)\big| |N^x_t| dt\Big] \le C(\e^n_1 + I(\e^{n-1}_1))\int_0^\infty \!\!\! e^{-\rho t}  \dbE[ |N^x_t|] dt\\
  &&\!\!\!\!\!\!\!\le C(\e^n_1 + I(\e^{n-1}_1))\Big[\int_0^\infty {1\over \sqrt{t}} e^{-\rho t+ C_1t}  dt\Big] = {C_2\over \sqrt{\rho- C_1}}(\e^n_1 + I(\e^{n-1}_1)),
\eeaa
where $C_1, C_2>0$ are generic constants independent of $\rho$, and we assume $\rho>C_1$. Set 
\bea
\label{beta0}
\rho_0 := C_1 + 9 C_2^2\q\mbox{so that} \q {C_2\over \sqrt{\rho_0- C_1}} = {1\over 3}.
\eea
Since $x$ is arbitrary, then for $\rho \ge \rho_0$ we obtain
\bea
\label{en1recursive}
\qq\q \e^n_1 \le {C_2\over \sqrt{\rho- C_1}}(\e^n_1 +  I(\e^{n-1}_1)) \le {1\over 3}(\e^n_1 +  I(\e^{n-1}_1)), \q\mbox{thus}\q \e^n_1 \le {1\over 2}  I(\e^{n-1}_1).
\eea
Moreover, by \reff{vnrep} and \reff{Hproperty} we have
\beaa
| v_x(x)| \le C(1+\| v_x\|_0)\int_0^\infty e^{-\rho t}\dbE[|N^x_t|]dt \le {C_2\over \sqrt{\rho - C_1}}(1+\| v_x\|_0) \le {1\over 3}(1+\| v_x\|_0).
\eeaa
This implies $\| v_x\|_0\le {1\over 2}$. Note further that $ v_x^0 \equiv 0$. Then, by \reff{en1recursive} and applying Lemma \reff{lem-superexponential} with $\eta_0=0$,  there exist $0<\eta<1$ and $C>0$ such that
\bea
\label{vxnrate}
 \e^n_1 \le C\eta^{2^n},\q \mbox{and thus}\q \| v_x^n\|_0 \le \| v_x\|_0 + \e^n_1 \le C. 
\eea

Next, it follows from  \reff{vnrep} that
$
|\D v^n(x)| \le \dbE\big[\int_0^\infty e^{-\rho t} \big|\D f^n(X^x_t)\big| dt\big].
$
Then, by \reff{fnfest} we have
\beaa
 |\D v^n(x)| \le C(\e^n_1 + I(\e^{n-1}_1))\int_0^\infty e^{-\rho t}  dt = {C\over \rho}(\e^n_1 + I(\e^{n-1}_1)).
\eeaa
Plug \reff{vxnrate} into it and note again that $\rho\ge \rho_0$, we obtain
\bea
\label{vnrate}
\| \D v^n\|_0  \le C\eta^{2^n}.
\eea
Moreover, similarly to \reff{pirate} we derive from \reff{vxnrate} the desired estimate for $\D \pi^n$.
 
 {\bf Step 3.} We now estimate  $\e^n_2 := \|  \D v_{xx}^n\|_0$. By \reff{vxnrep} we have, for any $x\in\hR^d$, 
\beaa
| \D v_{xx}^n(x)| \le \dbE\Big[\int_0^\infty e^{-\rho t}\big[ \big|\D f_x^n(X^x_t)\big| |\td X^x_t|| N^x_t|  + \big|\D f^n(X^x_t)\big||\td N^x_t|\big]   dt\Big].
\eeaa
Similarly to \reff{paxfn-f} we have, recalling the notation $\hat \e^n_i := \e^n_i + \e^{n-1}_i$, $i=1,2$, 
\beaa
 \Big|\D f_x^n(x)\Big|  \le C\big[1+\hat\e^n_1\big]\hat\e^n_2 + C\hat\e^n_1.
\eeaa
Thus, for possibly larger $C_1$ and $C_2$, by \reff{fnfest}, \reff{vxnrate}, and \reff{Ntxest} we have
\beaa
&&\dis | \D v_{xx}^n(x)| \le C\big[\hat\e^n_2 + \eta^{2^n} \big]\int_0^\infty e^{-\rho t}\dbE\big[|\td X^x_t|| N^x_t| +|\td N^x_t|\big]   dt\\
&&\dis \le C\big[\hat\e^n_2  +  \eta^{2^n}\big]\int_0^\infty {1\over \sqrt{t}}e^{-\rho t +C_1 t}dt \le {C_2\over \sqrt{\rho-C_1}} \big[\hat\e^n_2  +  \eta^{2^n}\big] \le {1\over 3}\big[\hat\e^n_2  +  \eta^{2^n}\big],
\eeaa
for all $\rho \ge \rho_0$, where $\rho_0$ is defined by \reff{beta0}, but with the $C_1, C_2$ here. Since $x$ is arbitrary, we have $\e^n_2\le  {1\over 3}\big[\hat\e^n_2  +  \eta^{2^n}\big] \le  {1\over 3}\big[\e^n_2 + \e^{n-1}_2 +  \eta^{2^n}\big]$.
 This implies 
\bea
\label{en2recursive}
  \e^n_2\le  {1\over 2}\e^{n-1}_2  + C\eta^{2^n},\q\mbox{and thus}\q \e^n_2\le  {\e^0_2\over 2^n}  + C\eta^{2^n}.
 \eea
 
 Moreover, by \reff{Hproperty} and noting from Step 2 that $\| v_x\|_0\le C$ for $\rho\ge \rho_0$, we have
 \beaa
 |f| = |H^v|\le C[1+| v_x|]\le C,\q |f_x| \le |H^v_x| + |H^v_z|| v_{xx}|\le C[1+\| v_{xx}\|_0]. 
 \eeaa
 Thus by \reff{vxnrep} and \reff{Ntxest} we have, again for $\rho\ge \rho_0$,
 \beaa
 | v_{xx}(x)| &\le&  C[1+\| v_{xx}\|_0]\int_0^\infty e^{-\rho t}\dbE[|\td X^x_t|| N^x_t|]  + C\int_0^\infty e^{-\rho t}\dbE[| \td N^x_t|]   dt\\
 &\le& {C_2\over \sqrt{\rho-C_1}}[1+\| v_{xx}\|_0] \le {1\over 3} [1+\| v_{xx}\|_0].
 \eeaa
 By the arbitrariness of $x$, we have $\| v_{xx}\|_0\le {1\over 3} [1+\| v_{xx}\|_0]$ and thus $\| v_{xx}\|_0\le {1\over 2}$. Note further that  $v^0_{xx}  \equiv 0$. Then $\e^0_2 = \|v_{xx}\|_0 \le {1\over 2}$, and thus it follows from \reff{en2recursive} that $\e^n_2 \le {C\over 2^n}$. This, together with \reff{vxnrate} and \reff{vnrate}, proves \reff{conv2}.
 \qed

\subsection{Proof of Theorem \ref{thm-main1} (ii)}
\label{sect-small}
We now prove \reff{conv1} for arbitrary $\rho$. Let $\rho_1>\rho$ be a large constant which will be specified later. We remark that here we allow $\rho_1$ to depend on $\rho$.  We proceed in two steps. 

{\bf Step 1.} In this step we estimate $\|v^n\|_2$. Note that we may rewrite \reff{PDEn} as
\beaa
\rho_1 v^{n}= \frac{1}{2} \sigma\sigma^\top :   v_{xx}^{n} + H^{n-1}_z \cd [v_x^{n} -  v_x^{n-1} ] + H^{n-1} + (\rho_1 - \rho) v^n.
\eeaa
Denote
\beaa
L^n_1 := \| v_x^n\|_0,\q L^n_2 := \| v_{xx}^n\|_0,\q \tilde f^n := f^n + (\rho_1-\rho) v^n.
\eeaa
 First, similarly to \reff{vxnrep} we have 
\beaa
 v_x^{n}(x)=\dbE\Big[\int_0^\infty e^{-\rho_1 t}\tilde f^n(X^x_t)  N_t^xdt\Big].
\eeaa
By Proposition \ref{prop-vn} (i) we have $\|v^n\|_0\le {C\over \rho}$.  Then, by \reff{vnrep} and Lemma \ref{lem-H}, 
\bea
\label{tildefnest}
|\tilde f^n(x)|\le C(L^n_1 + L^{n-1}_1) + C_\rho(\rho_1-\rho).
\eea
Thus, by \reff{Ntxest},
\beaa
| v_x^{n}(x)|&\le& C\big[L^n_1 + L^{n-1}_1 + C_\rho(\rho_1-\rho) \big]\int_0^\infty e^{-\rho_1 t} \dbE[|N_t^x|]dt\\
&\le&C\big[L^n_1 + L^{n-1}_1 + C_\rho(\rho_1-\rho) \big] \int_0^\infty {1\over \sqrt{t}} e^{-\rho_1 t+Ct}dt \\
&\le& {C_2\over \sqrt{\rho_1-C_1}}(L^n_1 + L^{n-1}_1) + {C_\rho(\rho_1-\rho) \over \sqrt{\rho_1-C_1}}.
\eeaa
Here $C_1, C_2$ are generic constants independent of $\rho$. Now set $\rho_1>\rho$ large enough as in \reff{beta0} such that ${C_2\over \sqrt{\rho_1-C_1}}\le {1\over 3}$. Then, by the arbitrariness of $x$, we obtain
\beaa
L^n_1 \le  {1\over 3}(L^n_1 + L^{n-1}_1) + C_\rho,\q\mbox{and thus}\q L^n_1 \le {1\over 2} L^{n-1}_1 + C_\rho.
\eeaa
 Note further that $ v_x^0=0$ and thus $L^0_1=0$. Then by standard arguments we have
 \bea
 \label{Ln1est}
 L^n_1 \le 2^{-n}L^0_1   + C_\rho \le C_\rho.
 \eea
 
 Next, similarly to \reff{vxnrep} we have
 \bea
 \label{paxxvnrep}
\qq  v_{xx}^{n}(x) = \dbE\Big[\int_0^\infty e^{-\rho_1 t}\big[ N^x_t(\tilde f_x^n(X^x_t))^\top \td X^x_t    +\tilde f^n(X^x_t) \td N^x_t\big]   dt\Big].
 \eea
By \reff{tildefnest} and \reff{Ln1est} its is clear that $|\tilde f^n(x)|\le C_\rho(\rho_1 - \rho +1)$. Moreover, following similar arguments as in \reff{paxfn-f} and by using \reff{Ln1est} again, we have
\beaa
 |\tilde f_x^n (x)| \le | f_x^n (x)| + (\rho_1-\rho) | v_x^n|  \le C_\rho\big[L^n_2 + L^{n-1}_2 + \rho_1-\rho + 1\big].
\eeaa
Then, by \reff{Ntxest} and by the arbitrariness of $x$,
\beaa
L^n_2 &\le&C_\rho\big[L^n_2 + L^{n-1}_2 + \rho_1-\rho + 1\big] \dbE\Big[\int_0^\infty {1\over \sqrt{t}}e^{-\rho_1 t+Ct}dt\Big]\\
&\le& {C_\rho\over \sqrt{\rho_1-C_1}}\big[L^n_2 + L^{n-1}_2 + \rho_1-\rho + 1\big]\le {1\over 3}\big[L^n_2 + L^{n-1}_2 + \rho_1-\rho + 1\big],
\eeaa
where in the last inequality we set $\rho_1:= C_1 + 9|C_\rho|^2$  so that $ {C_\rho\over \sqrt{\rho_1-C_1}} = {1\over 3}$. Note that $L^0_2 = \| v_{xx}^0\|_0=0$. By standard arguments this implies that
\bea
\label{Ln2est}
L^n_2 \le {1\over 2}\big[ L^{n-1}_2 + \rho_1-\rho + 1\big],~\mbox{and thus},~ L^n_2 \le C(\rho_1-\rho + 1)\le C_\rho.
\eea
{\bf Step 2.} We now prove the desired convergence. First, by the monotonicity and boundedness of $v^n$, there exists bounded $v^*$ such that $v^n\uparrow v^*$. By \reff{Ln1est} $\{v^n\}_{n\ge 0}$ are equicontinuous, then the above convergence is uniform on compacts. Next, by \reff{Ln1est} and \reff{Ln2est} we see that $\{ v_x^n\}_{n\ge 0}$ are bounded and equicontinuous. Then by applying Arzela-Ascoli Theorem there exist a subsequence $\{n_k\}_{k\ge1}$ such that $ v_x^{n_k}$ converge uniformly on compacts. Note that differentiation is a closed operator, and since $v^n\to v^*$, we must have $ v_x^{n_k}\to  v_x^*$. This implies that the limit of the subsequence $\{ v_x^{n_k}\}_{k\ge 1}$ is unique, then we must have the convergence of the whole sequence $ v_x^n$, namely $ v_x^n \to  v_x^*$ uniformly on compacts. In particular, this implies that
\beaa
f^n\to f^* ~\mbox{uniformly on compacts, where}\q f^*(x) := H(x,  v_x^*).
\eeaa
 Moreover, by \reff{Ln2est} it is clear that $ v_x^*$, whence $f^*$, is uniformly Lipschitz continuous.
 
Note that $v^n$ is the classical solution of the PDE:
$
 \rho v^n = {1\over 2} \si\si^\top :  v_{xx}^n + f^n.
$
 Let $\tilde v$ denote the unique viscosity solution of the PDE:
 \bea
 \label{tildevviscosity}
 \rho \tilde v = {1\over 2} \si\si^\top : \pa_{xx} \tilde v + f^*.
 \eea
 By the stability of the viscosity solution, we see that $v^* = \lim_{n\to\infty} v^n$ is a viscosity solution of \reff{tildevviscosity}. Moreover, since $ v_x^*$ is (Lipschitz) continuous, for any smooth test function $\f$ of $v^*$ at $x$ in the definition of viscosity solution, we must have $ \f_x(x) =  v_x^*(x)$ and thus $H(x, \f_x(x)) = f^*(x)$, then $v^*$ is also a viscosity solution of the PDE:
 \beaa
 \rho \tilde v = {1\over 2} \si\si^\top : \pa_{xx} \tilde v + H(x, \pa_x \tilde v).
 \eeaa
 This PDE identifies with \reff{HJBv}, then by the uniqueness of its viscosity solution, we obtain $v^* = v$. That is, $(v^n,  v_x^n) \to (v,  v_x)$ uniformly on compacts.
 
It remains to prove the desired convergence of $ v_{xx}^n$.\footnote{Provided a uniform H\"{o}lder continuity of $ v_{xx}^n$, which could be obtained by PDE estimates, one may derive the convergence of $ v_{xx}^n$ by the same compactness argument at above.  Here we provide a probabilistic proof for the convergence directly, which does not require the H\"{o}lder continuity of $ v_{xx}^n$.} To this end, we shall first introduce another representation formula for $ v_{xx}^n$. Let us recall \reff{DXtx} and denote
\bea
\label{R}
R^x_t := N^x_t(N^x_t)^\top -{1\over t} \int_0^t D_s N^x_t \check \si(X^x_s) \td X_s ds + \td N^x_t,
\eea
where $D_s N^x_t$ is the Malliavin derivative, see \cite{Nualart}, and $\check\si := \si^{-1}$. Note that, denoting by $D_s^i N^x$ (resp. $\td_i X^x$) the $i$-th column of $D_s N^x$ (resp. $\td X^x$),
\beaa
&\dis D^i_s N^x_t = {1\over t}\check\si(X^x_s) \td_i X^x_s + {1\over t} \int_s^t \Big(\check \si_{x_j}(X^x_l) D^i_s X^{x,j}_l \td X_s + \check\si(X^x_s) D^i_s \td X_l\Big)^\top dW_l,\\ 
&\dis D^i_s X^x_t = \si^i(X^x_s) + \int_s^t \si_{x_j}(X^x_l) D^i_s X^{x,j}_ldW_l,\\
&\dis D^i_s \td X^x_t = \si^i_x(X^x_s) \td X^x_s +  \int_s^t \Big(\si^j_{x_kx}(X^x_l) D^i_s X^{x,k}_l  \td X^x_l +  \si^j_x(X^x_l) D^i_s \td X^x_l\Big) dW^i_l.
\eeaa
Here we used the Einstein summation again. Fix $s$, and consider $D^i_s X^x, D^i_s \td X^x$ as the solution to the above linear SDE systems for $t\in [s, \infty)$. One can easily 
 check that 
 \bea
 \label{Rest}
 \dbE[|D_s N^x_t|^4]\le {C\over t^4} e^{Ct}, \q \mbox{and thus}\q \dbE\big[|R^x_t|^2] \le {C\over t^2} e^{Ct}.
 \eea
  Then, for any $\phi\in C^0_b(\dbR^d; \dbR)$, by \cite[Chapter 2]{Zhangthesis} we have\footnote{As in Footnote \ref{footnote-rep}, when $d=1$ and $\si \equiv 1$, we have $\td X^x_t = 1$, $N^x_t = {W_t\over t}$, $D_s N^x_t = {1\over t}$, $\td N^x_t =0$, then $R^x_t = {W_t^2 - t  \over t^2}$, and thus
\beaa
 \phi_{xx}(X^x_t) = \pa_{xx} \int_\dbR \phi(y) {1\over \sqrt{2\pi t}} e^{-{(y-x)^2\over 2t}}dy =  \int_\dbR \phi(y) {1\over \sqrt{2\pi t}} e^{-{(y-x)^2\over 2t}}{(y-x)^2 - t\over t^2}dy = \dbE[\phi(X^x_t) R^x_t]. 
\eeaa
 }
 \bea
 \label{Rrep}
 \pa_{xx} \hE[\phi(X^x_t)] = \hE\big[\phi(X^x_t)R^x_t\big],
 \eea
 Thus, for any $\d>0$ small, we may rewrite  \reff{paxxvnrep} as
 \beaa
   v_{xx}^{n}(x) = \dbE\Big[\int_0^\d \!\! e^{-\rho_1 t}\big[N^x_t (\tilde f_x^n(X^x_t))^\top\td X^x_t    +\tilde f^n(X^x_t) \td N^x_t\big]   dt+\!\! \int_\d^\infty \!\!\! e^{-\rho_1 t}\tilde f^n(X^x_t) R^x_t   dt\Big].
  \eeaa
  We remark that $\dbE[|\tilde f^n(X^x_t) R^x_t|] \le {C\over t}e^{Ct}$ which is not integrable around $t=0$, so at above we use different representations for small $t$ and large $t$.  Similarly we have
  \beaa
 &\dis  v_{xx}(x) = \dbE\Big[\int_0^\d \!\! e^{-\rho_1 t}\big[N^x_t (\tilde f_x(X^x_t))^\top\td X^x_t      +\tilde f(X^x_t) \td N^x_t\big]   dt+\!\!\int_\d^\infty \!\!\! e^{-\rho_1 t}\tilde f(X^x_t) R^x_t   dt\Big];\\
&\dis \mbox{where}\q \tilde f(x) = H^v(x) + (\rho_1-\rho) v.
 \eeaa

Now fix a compact set $K \Subset \hR^d$, and assume that $K\subseteq B_{M_0}(0)$ for some $M_0>0$.  
 Denote, for all $M>M_0$,
 \beaa
 \e^n_M  := \sup_{|x|\le M} \big[|\D v^n(x)| + |\D v^n_x(x)|\big]\to 0,\q\mbox{as}~ n\to \infty.
 \eeaa
Then, by \reff{Ntxest} and \reff{Rest}, one can easily see that, for any $x\in K$,
 \beaa
 &&\dis | \D v_{xx}^n(x)| \le  \hE\Big[\int_\d^\infty e^{-\rho_1 t}\big|\D\tilde f^n(X^x_t)\big| |R^x_t|   dt \\ 
&&\dis \qq+\int_0^\d e^{-\rho_1 t}\big[\big( |\tilde f_x^n|+|\tilde f_x|\big)(X^x_t) |\td X^x_t|| N^x_t|  +\big(|\tilde f^n|+|\tilde f|\big)(X^x_t) |\td N^x_t| \big]dt \Big]\\
&&\dis \le  C_\rho \int_0^\d\!\! e^{-\rho_1 t}\dbE\big[|\td X^x_t || N^x_t| + |\td N^x_t| \big]dt  \\
&&\dis \qq+ \int_\d^\infty\!\! e^{-\rho_1 t}\dbE\big[\big(\e_M^n  + C_\rho\1_{\{|X^x_t|\ge M\}} \big)|R^x_t| \big]  dt\\
&& \dis  \le C_\rho \sqrt{\d} + C \e_M^n \int_\d^\infty {1\over t}e^{-\rho_1 t+C t}  dt +  {C_\rho\over M} \dbE\Big[\int_\d^\infty e^{-\rho_1 t}|X^x_t| |R^x_t|   dt\Big]\\
&&\dis\le  C_\rho \sqrt{\d} +C_\rho \e_M^n  \ln {1\over \d}  +  {C_\rho\over M} (1+|x|)\int_\d^\infty {1\over t}  e^{-\rho_1 t+Ct} dt\\
&&\dis\le  C_\rho \sqrt{\d} +C_\rho \e_M^n  \ln {1\over \d}  +  {C_\rho\over M} (1+|M_0|) \ln {1\over \d}.
 \eeaa
 Thus
 \beaa
 \sup_{x\in K} | \D v_{xx}^n(x)| \le C_\rho \sqrt{\d} +C_\rho \e_M^n  \ln {1\over \d}  +  {C_\rho\over M} (1+|M_0|) \ln {1\over \d}.
 \eeaa
 Fix $M$, $\d$ and send $n\to \infty$, we obtain
 \beaa
 \limsup_{n\to\infty}  \sup_{x\in K} | \D v_{xx}^n(x)| \le C_\rho \sqrt{\d}   +  {C_\rho\over M} (1+|M_0|) \ln {1\over \d}.
\eeaa
By first sending $M\to \infty$ and then $\d\to 0$, we obtain the desired estimate: 
\beaa
\limsup_{n\to\infty}  \sup_{x\in K} | \D v_{xx}^n(x)|=0.
\eeaa

\vspace{-8mm}
\qed

\ms
\begin{rem}
\label{rem-scalar}
When $d=1$, the uniform estimate for $\| v_{xx}^n\|_0$ in Step 1 and the convergence of $ v_{xx}^n$ in  Step 2 become trivial. Indeed, in this case we have
\beaa
 v_{xx}^n = {2\over \si^2}\big[\rho v^n - H^{n-1}_z[ v_x^n- v_x^{n-1}] - H^{n-1}\big].
\eeaa
Then the desired boundedness and convergence of $ v_{xx}^n$ follow directly from those of $v^n, v_x^n$. In particular, when $\rho\ge \rho_0$, by \reff{conv2} we also have $\|\D v^n\|_2 \le C\eta^{2^n}$. We shall use this feature to study a diffusion control model in the next section.
\end{rem}

\begin{rem}
\label{rem-diverge}
The convergence in this case relies heavily on the fact that $v^n$ is monotone and hence converges. When $\rho$ is small, in general Picard iteration may not converge, not to mention rate of convergence, as we see in the following example.  
\end{rem}

\begin{eg}
\label{eg-diverge}
Let $d=1$. Consider the following (linear) PDE with unique bounded classical solution $v\equiv 0$:
\beaa
 \rho v = {1\over 2}  v_{xx} +  v_x.
 \eeaa
 Set $v^0(x) := -\cos x$, and define $v^n$ recursively by:
$
 \rho v^n = {1\over 2}  v_{xx}^n +  v_x^{n-1}.
$
 Then 
 \bea
 \label{paxvn0}
    v_x^n(0) = \left\{\ba{lll} 0,\qq \mbox{$n$ is even}; \ss\\ (-1)^{n-1\over 2}(\rho + {1\over 2})^{-n}, \q \mbox{$n$ is odd}. \ea\right.
  \eea
In particular,  $| v_x^{2m+1}(0)| \to \infty$ as $m\to \infty$, whenever $\rho <{1\over 2}$.
 \end{eg}
 \proof By \reff{vxnrep}, we have
 \beaa
&\dis  v_x^n(x) = \dbE\Big[\int_0^\infty e^{-\rho t_1}  v_x^{n-1}(x+ W_{t_1}){W_{t_1}\over t_1}dt_1\Big],\\
&\dis  v_x^{n-1}(x+W_{t_1}) =  \dbE\Big[\int_0^\infty e^{-\rho t_2}  v_x^{n-2}(x+ W_{t_1+t_2}){W_{t_1+t_2}-W_{t_1}\over t_2}dt_2\Big|\cF^W_{t_1}\Big].
 \eeaa 
 Plug the second formula into the first one, we get
  \beaa
&\dis  v_x^n(x) = \dbE\Big[\int_{\dbR_+^2} e^{-\rho (t_1+t_2)}  v_x^{n-2}(x+ W_{t_1+t_2}){W_{t_1}\over t_1}{W_{t_1+t_2}-W_{t_1}\over t_2}dt_2dt_1\Big].
 \eeaa 
 Repeat the arguments and note that $ v_x^0(x) = \sin x$, we obtain
 \beaa
 v_x^n(x) = \dbE\Big[\int_{\dbR_+^n} e^{-\rho T_n} \sin(x+ W_{T_n}){W_{t_1}\over t_1}\cds{W_{T_n}-W_{T_{n-1}}\over t_n}dt_n\cds dt_1\Big],
 \eeaa 
 where $T_i := t_1+\cds+t_i$. 
Note that $E[\cos W_t {W_t \over t}]=0$, and 
 \beaa
 \dbE[\sin W_t{W_t \over t}] = \sum_{k=0}^\infty (-1)^k \dbE\big[{W_t^{2k+1}\over (2k+1)!}{W_t\over t} \big] = \sum_{k=0}^\infty (-t)^k {(2k+1)!!\over (2k+1)!}  =  \sum_{k=0}^\infty {(-t)^k \over 2^k k!} = e^{-{t\over 2}}.
 \eeaa 
Let Im denote the imaginary part of a complex number. Then
\beaa
&&\dis  v_x^n(0) = {\rm Im}~ \dbE\Big[\int_{\dbR_+^n} e^{-\rho T_n } e^{\sqrt{-1} W_{T_n}}{W_{t_1}\over t_1}\cds{W_{T_n}-W_{T_{n-1}}\over t_n}dt_n\cds dt_1\Big]\\
&&\dis= {\rm Im} \int_{\dbR_+^n}  \prod_{i=1}^n e^{-\rho t_i }\dbE\Big[ e^{\sqrt{-1} (W_{T_i}-W_{T_{i-1}})}{W_{T_i}-W_{T_{i-1}}\over t_i}\Big] dt_n\cds dt_1\\
&&\dis= {\rm Im} \int_{\dbR_+^n}  \prod_{i=1}^n e^{-\rho t_i }\dbE\big[ e^{\sqrt{-1} W_{t_i}}{W_{t_i}\over t_i}\big] dt_n\cds dt_1 = {\rm Im} \Big(\int_0^\infty  e^{-\rho t }\dbE\Big[ e^{\sqrt{-1} W_t}{W_t\over t}\Big] dt\Big)^n\\
&&\dis= {\rm Im}\Big( \int_0^\infty  e^{-\rho t } \sqrt{-1} e^{-{t\over 2}} dt\Big)^n = {\rm Im}\Big({\sqrt{-1}\over \rho + {1\over 2}} \Big)^n.
 \eeaa 
 This implies \reff{paxvn0} immediately.
 \qed

\section{The Scalar Case with Diffusion Control}
\label{sect-diffusion}
\setcounter{equation}{0}

In this section we consider the diffusion control  case, where the corresponding HJB equation becomes  fully nonlinear. It has been widely recognized that the general case is much more challenging than the drift control case, which we shall leave to future research. In this section we consider only the {\it one-dimensional case}, i.e,  $d=1$.  Recall Remark \ref{rem-scalar}.

Consider the setting in \S\ref{sect-infinite}, our entropy-regularized problem is:
\bea
\label{scalarV}
\left.\ba{lll}
\dis dX^\pi_t = \tilde b(X^\pi_t, \pi(t,X^\pi_t))dt + \sqrt{\widetilde {\si^2}(X^\pi_t, \pi(t, X^\pi_t))} dW_t;\ss\\
\dis J(x, \pi):= \hE\Big[\int_0^\infty e^{-\rho t} \big[\tilde r(X^\pi_t, \pi(s, X^\pi_s)) + \l \cH(\pi(t, X^\pi_t))\big]dt\Big|X^\pi_0=x\Big], \ss\\
\dis v(x) := \sup_{\pi\in \cA_{Lip}} J(x, \pi).
\ea\right.
\eea
Here in the above, to simplify the arguments we restrict the admissible relaxed controls to $\cA_{Lip}$, the set of 
 $\pi\in \cA$ that is Lipschitz continuous in $x$, so that the $X^\pi$ above has a unique strong solution. Furthermore, in this section we shall assume:
 
 \ss
\begin{assum}
\label{assum-infinite2} $d=1$; $b, \si, r$ satisfy  Assumption \ref{assum-finite}, with $\si$ depending on $a$; and  $b, \si$ are uniformly continuous in $a$, uniformly in $x$.
\end{assum}

It is well-known that, in this case  $v$ satisfies the  fully-nonlinear HJB equation:
\bea
\label{scalarHJBv}
\left.\ba{c}
\dis \rho v = H(x,  v_x,  v_{xx}), \q x\in\hR,\q\mbox{where}\ss\\
\dis
H(x, z, q) := \sup_{\pi\in \cP_0(A)} \Big[{1\over 2} \widetilde {\si^2}(x, \pi) q  + \tilde b(x, \pi) z+ \tilde r(x, \pi) + \l \cH(\pi)\Big].
\ea\right.
\eea
Moreover, the maximizer of the Hamiltonian $H$ has the Gibbs form $\G$:
\bea
\label{scalarHGamma}
\left.\ba{c}
\dis \pi^*(x,a) = \G(x, v_x, v_{xx}, a),\q \G(x, z, q, a) := \frac{\g(x,z, q, a)}{\int_A \g(x,z, q, a') d a'},\ss\\
\dis \mbox{where}~ 
\g(x,z, q, a):= \exp \Big(\frac{1}{\lambda}[{1\over 2}\si^2(x,a)q+b(x,a)z+r(x, a)]\Big),
\ea\right.
\eea 
and consequently, we have
\bea
\label{scalarH}
H(x,z, q)=\lambda \ln \Big(\int_A \g(x,z, q, a) d a\Big).
\eea
We have the following simple extension of Lemma \ref{lem-H}: 
\begin{lem}
\label{lem-scalarH}
Let Assumption \ref{assum-infinite2} hold. Then $H$ is twice continuously differentiable in $(x, z,q)$;  jointly convex in $(z,q)$; and, for some  constant $C>0$,
\bea
\label{scalarHproperty}
\qq |H_z|, |H_q|, |H_{zz}|,  |H_{zq}|, |H_{qq}|\le C, ~  H_q \ge {1\over C},~|H(x,z,q)|\le C[1+|z|+|q|].  
\eea
\end{lem}

From Lemma \ref{lem-scalarH} we see that $H$ is convex, strictly increasing, and also has linear growth in $q$. The following observation about the asymptotic behavior of $H$ in $q$ is crucial for our convergence analysis on the recursive PDEs.

\begin{prop}
\label{prop-scalarH}
Assume that Assumption \ref{assum-infinite2} is in force. Then, for any $\e>0$, there exists $ C_\e>0$ such that
\bea
\label{HpaqH}
 |h(x, z, q)|\le  \e |q| + C|z|+ C_\e,\q\mbox{where}\q  h:= H-H_z z -  H_q q. 
  \eea
 In particular, if $b, \si$ are uniformly H\"{o}lder continuous in $a$, uniformly in $x$, then  
\bea
\label{HpaqHln}
 |h(x, z, q) |\le  C \big[1+|z| + \ln (1+|q|)\big].
  \eea
\end{prop}
\proof  We prove the result only for $q>0$. The case $q<0$ can be proved similarly\footnote{Actually this case is not needed, because later on we can easily show that $ v_{xx}^n\ge -C$ for all $n$.}.

First, since $H$ is jointly convex in $(z,q)$, we have 
\bea
\label{HpaqH1}
 h(x, z, q)  \le H(x, 0, 0) \le C.
 \eea
 To see the opposite inequality, denote $\th(x,z,q,a) := {1\over 2}\si^2(x,a)  + b(x,a){z\over q}$. When there is no confusion, we omit the variables $(x,z,q)$ in $\th$, $\g$.  Then, by (\ref{scalarHGamma}), \reff{scalarH} we have 
 \bea
 \label{HpaqH2}
-h(x, z, q)= {q \int_A \th(a) \g( a) da\over  \int_A \g(a) da} - \lambda \ln \Big(\int_A \g(a) d a\Big).
  \eea
 Denote $C_1 := \sup_{a\in A, x\in \dbR} {2|b(x,a)|\over |\si^2(x,a)|}<\infty$. When $q \le 2C_1|z|$, the result is obviously true. We now assume $q \ge 2C_1|z|$, which implies $\th(a) \ge {1\over C_2}>0$, for some $C_2>0$,  for all $(x, a)$.  Fix $(x, q, z)$ and denote 
 \beaa
 \ol\th := \sup_{a\in A}\th(a) > {1\over C_2}, \q A_\e := \big\{a\in A: \th(a) \ge  \ol\th -\e\big\}.
 \eeaa
  Let $a^\e\in A$ be such that $\th(a^\e) \ge \ol \th-{\e\over 2}$. Since $b, \si$ are uniformly continuous in $a$,  uniformly in $x$, and $\si$ is bounded and  $q \ge 2C_1|z|$, we see that $\th$ is also uniformly continuous in $a$.  Thus there exists $\d_\e>0$, independent of $(x,z,q)$, such that $\th(a) \ge \ol\th -\e$ whenever $|a-a^\e|\le \d_\e$. That is, $A_\e \supset A \cap B_{\d_\e}(a^\e)$. Then, since $A$ has smooth boundary, there exists $\mu_\e>0$, depending only on the model parameters,  such that $|A_\e|\ge \mu_\e$ for all $(x, z, q)$ with $q \ge 2C_1|z|$.   Note that
 \beaa
 \left.\ba{lll}
  \int_{A}  \g( a) da \neg\neg& \neg\ge\neg &\neg\neg \int_{A_\e}  \g( a) da \ge e^{ {(\ol\th-\e) q\over \l} - C}\mu_\e; \ss\\
 {\int_{A\backslash A_\e}\g(a) da\over \int_{A_\e} \g(a) da}  \neg\neg&\neg \le \neg&\neg\neg {\int_{A\backslash A_\e} \g(a) da\over \int_{A_{\e\over 2}} \g(a) da}  \le {e^{{(\ol\th - \e) q\over \l }+ C}|A\backslash A_\e| \over e^{{(\ol\th -{\e\over 2})q\over \l} - C}|A_{\e\over 2}|}  \le e^{-{\e \over 2\l} q + C}{|A|  \over \mu_{\e\over 2}};\ss\\
 \int_A\th(a)  \g(a) da   \neg\neg&\neg \le \neg&\neg\neg  \ol \th \int_{A_\e}  \g(a) da  +  (\ol \th -\e ) \int_{A\backslash A_\e}  \g(a) da\ms \\
\neg\neg&\neg \le \neg&\neg\neg 
 \big[\ol \th + (\ol \th -\e )e^{-{\e \over 2\l} q + C}{|A|  \over \mu_{\e\over 2}}\big] \int_{A_\e} \g(a) da.
 \ea\right.
  \eeaa
  Then, by \reff{HpaqH2}  and assuming without loss of generality that  $\e<{1\over C_2}$ so that $\ol\th -\e>0$,
 \bea
 \label{paqHHest}
 \left.\ba{lll}
  \dis -h(x, z, q)
\le q\Big[\ol \th + (\ol \th -\e)e^{-{\e \over 2\l} q + C}{|A|  \over \mu_{\e\over 2}}\Big]  -\l \ln \big(e^{ {(\ol\th-\e) q\over \l} - C}\mu_\e\big)\\
\dis =q\Big[\ol \th + (\ol \th -\e )e^{-{\e \over 2\l} q + C}{|A|  \over \mu_{\e\over 2}}\Big] -\big[(\ol\th-\e) q - \l C + \l \ln \mu_\e\big]\le \e q + C_\e.
\ea\right.
 \eea
 This, together with \reff{HpaqH1}, proves  \reff{HpaqH}.
 
 Finally, if $b, \si$ are H\"{o}lder-$\beta$ continuous in $a$, then we can easily see that $\mu_\e \ge {1\over C} \e^{1\over \beta}$. Thus, from the second line of \reff{paqHHest} we have
 \beaa
 -h(x, z, q) 
&\le &   q\Big[\ol \th + (\ol \th -\e )e^{-{\e \over 2\l} q + C}C  \e^{-{1\over \beta}}\Big] -\Big[(\ol\th-\e) q - C + \l \ln \e^{1\over \beta}\Big]\\
&\le& \e q + C q e^{-{\e \over 2\l} q}\e^{-{1\over \beta}} - C\ln \e +  C.
 \eeaa
 Set $\e := 2\l (1+{1\over \beta}) {\ln q\over q}$. Then, assuming $q\ge e$  without loss of generality, 
 \beaa
-h(x, z, q)
&\le& C\ln q + C q e^{-(1+{1\over \beta})\ln q} ({q\over \ln q})^{1\over \beta} + C\ln q - C \ln\ln q +  C\\
& =& C\ln q + C (\ln q)^{-{1\over \beta}} + C\ln q - C \ln\ln q +  C \le C\ln q + C.
 \eeaa
This, together with  \reff{HpaqH1} and the estimate in the case $q \le 2C_1|z|$, leads to \reff{HpaqHln}.
 \qed

\ms
For the PIA corresponding to \reff{scalarHJBv}, we set $v^0$ the same as in \S\ref{sect-infinite}, and for $n\ge 1$, define $\pi^{n}(x, a):=  \Gamma\left(x,v_x^{n-1}(x) , v_{xx}^{n-1}(x), a\right)$ and $v^{n}(x):=J(x, \pi^n)$. Then,  using (\ref{scalarHGamma}) and (\ref{scalarH}), it is easy to check that $v^{n}$ satisfies the following recursive linear PDE: 
\bea 
\label{scalarPDEn}
\left.\ba{lll}
\dis \rho v^{n} =   H^{n-1}_q[v_{xx}^{n}-  v_{xx}^{n-1}]  + H^{n-1}_z [ v_x^{n}-  v_x^{n-1} ] +H^{n-1}\ms\\
\dis\qq =   H^{n-1}_q v_{xx}^{n}+ H^{n-1}_z v_x^{n}+h^{n-1},\ms\\
\dis\mbox{where}\q \f^n(x) := \f(x, v^n_x, v^n_{xx})~\mbox{for}~ \f = H_q, H_z, h, \mbox{etc.}
\ea\right.
\eea
Here we abuse the notation $\phi^n$ with \reff{unrep}, \reff{vxnrep}. As in Proposition \ref{prop-vn}, we have

\ss
\begin{prop}
\label{prop-scalarvn}
Let Assumption \ref{assum-infinite2} hold. Then

(i) For each $n\ge 1$, $v^{n}\in C^{2}_b(\dbR^d; \dbR)$ is a classical solution of \reff{scalarPDEn}, with $v^n_{xx}$ being uniformly Lipschitz continuous, and hence $\pi^n\in \cA_{Lip}$.

(ii) $v^{n}$ is increasing in $n$ and $v^{n} \le {1\over \rho}\big[ C_0 +\l (\ln |A|)^+\big]$.
\end{prop}
\proof (ii) is standard and we shall only prove (i). Note that $v^0$ is a constant and thus satisfies all the requirements. Assume $v^{n-1}$ satisfies the requirements in (i). Then by the standard Schauder's estimate (cf. \cite{Krylov}), one can easily see that  $v^{n}\in C^2_b(\dbR^d; \dbR)$. Moreover, since $d=1$, by  \reff{scalarPDEn} we may write down $ v_{xx}^n$ explicitly:
\bea
\label{scalarPDEn1}
\left.\ba{c}
 v_{xx}^{n} = {1\over   H^{n-1}_q}\big[ \rho v^{n}   - H^{n-1}_z  v_x^{n}-h^{n-1} \big].
\ea\right.
 \eea
Then by \reff{scalarHproperty} it is clear that $v^n_{xx}$ is uniformly Lipschitz continuous.
\qed

\ss
 Our main result is as follows.

\begin{thm}
\label{thm-main2}
Under Assumption \ref{assum-infinite2},  $v^n\to v$ in $C^2$ and $\pi^n \to \pi^*$, uniformly on compacts in the sense of \reff{conv1}. Moreover, $v\in C^3_b(\dbR;\dbR)$ and  thus $\pi^*\in \cA_{Lip}$.
\end{thm}

\ms

\proof Again, we proceed in several steps. Denote
\beaa
L^n_1 := \| v_x^n\|_0, \q L^n_2:= \| v_{xx}^n\|_0,\q \ol L^n_1 := \sum_{k=1}^n {L^k_1\over 3^{n-k+1}}.
\eeaa

{\bf Step 1.} 
 For any $\e>0$, by \reff{scalarHproperty} and \reff{HpaqH}, it follows from \reff{scalarPDEn1} that
   \beaa
 L^n_2 \le {\e\over C_1} L^{n-1}_2  + C(L^n_1 + L^{n-1}_1) +C_{\rho,\e} \le {1\over 3} L^{n-1}_2  + C(L^n_1 + L^{n-1}_1) +C_\rho,
 \eeaa
 where we set $\e := {C_1\over 3}$ in the second inequality. Then by standard arguments we have
  \bea
 \label{scalar-paxxvnest}
L^n_2\le  C\ol L^n_1 +C_\rho.
 \eea

 {\bf Step 2.} Let $\rho_1>0$ be a large constant. Rewrite \reff{scalarPDEn} as:
 \bea
 \label{scalar-rho1}
\!\qq~~ \rho_1 v^n \!= {1\over 2}  v_{xx}^n +  f_n + (\rho_1-\rho) v^n,  ~ \! f_n := H^{n-1}_qv_{xx}^{n}  + H^{n-1}_z v_x^n + h^{n-1} \! -{1\over 2}  v_{xx}^n.
 \eea
 Then, denoting $X^x_t:= x+W_t$, by \reff{rep-d=1} we have
 \bea
  \label{scalar-rho1vnx}
  v_x^n(x) := \dbE\Big[\int_0^\infty e^{-\rho_1 t} \big[f_n(X^x_t) + (\rho_1-\rho) v^n(X^x_t)\big] {W_t\over t}dt\Big].
 \eea
 By \reff{scalarHproperty}, \reff{HpaqH}, \reff{scalar-paxxvnest} and noting that $\Bar{L}_1^n={\Bar{L}_1^{n-1}+L_1^n\over 3}$ and $L^{n-1}_1 \le 3 \Bar{L}_1^{n-1}$, we have
 \beaa
 | v_x^n(x)| &\le& C\Big[L^n_2+ L^{n-1}_2+  L^{n}_1+ L^{n-1}_1+1 + C_\rho|\rho_1-\rho|\Big] \int_0^\infty e^{-\rho_1 t}{1\over \sqrt{t}}dt\\
&\le& {C\over \sqrt{\rho_1}}\Big[ \ol L^n_1 + \ol L^{n-1}_1+L^n_1+ L^{n-1}_1+C_\rho + C_\rho|\rho_1-\rho|\Big] \\
 &\le&{C_1\over \sqrt{\rho_1}}\Big[ L^n_1+\ol L^{n-1}_1\Big] + {C_\rho\over \sqrt{\rho_1}}\big[|\rho_1-\rho|+1\big].
 \eeaa
Setting $\rho_1 = 16C_1^2$ and by the arbitrariness of $x$, we get 
\beaa
L^n_1 \le {1\over 4} (L^n_1+\ol L^{n-1}_1) + C_\rho,\q\mbox{and thus}\q L^n_1 \le {1\over 3} \ol L^{n-1}_1 + C_\rho
\eeaa
Note that 
\beaa
\ol L^n_1 = {1\over 3} L^n_1 + {1\over 3} \ol L^{n-1}_1 \le {1\over 9} \ol L^{n-1}_1 + {1\over 3} \ol L^{n-1}_1 + C_\rho\le {1\over 2} \ol L^{n-1}_1 + C_\rho.
\eeaa
 This, together with \reff{scalar-paxxvnest}, implies immediately that 
\bea
\label{scalar-paxvnest}
\ol L^n_1 \le C_\rho,\q\mbox{and thus}\q L^n_1\le C_\rho,\q L^n_2  \le C_\rho.
\eea

{\bf Step 3.} Follow the arguments in the beginning of Step 2 in \S\ref{sect-small}, we have $(v^n,  v_x^n)\to (v^*,  v_x^*)$ uniformly on compacts for some function $v^*\in C^1_b(\dbR; \dbR)$ such that $ v_x^*$ is Lipschitz continuous. Fix an arbitrary $x_0$, and denote $q_*:= \liminf_{n\to\infty} v_{xx}^n(x_0) =  \lim_{k\to\infty} v_{xx}^{n_k}(x_0)$ for some subsequence $\{n_k\}_{k\ge 1}$, which may depend on $x_0$.  By \reff{scalarPDEn} we have, at $x_0$,
\bea
\label{scalarPDEn2}
  v_{xx}^{n} &=&  v_{xx}^{n-1} + {1\over  H^{n-1}_q} \Big[ \rho v^{n} - H^{n-1}_z  [ v_x^{n}-  v_x^{n-1} ] - H^{n-1}(x_0)\Big]\nonumber\\
&\le&  v_{xx}^{n-1} + {\rho v^* - H(x_0, v_x^*,   v_{xx}^{n-1})\over  H_q(x_0, v_x^*,  v_{xx}^{n-1})} + C\e_n,
\eea
where
\beaa
\e'_n:=\Big[|v^n-v^{*}|+| v_x^n -  v_x^{*}|\Big](x_0),\q  \e_n:= \e'_n+\e'_{n-1} \to  0,~\mbox{as}~n\to \infty.
\eeaa
Since $H$ is convex in $(z, q)$, by \reff{scalarPDEn} again we have
\bea
\label{scalarPDEn3}
\qq  \rho v^{n} = H^n  - \Big[ H^n - H^{n-1}- H^{n-1}_q[v_{xx}^{n}-  v_{xx}^{n-1}]  - H^{n-1}_z  [ v_x^{n}-  v_x^{n-1}]\Big]
\le H^n. 
\eea
Set $n=n_k+1$ and send $k\to \infty$, we have $\rho v^* \le   H(x_0, v_x^*,  q_*)$. Then, by \reff{scalarPDEn2},
\beaa
 v_{xx}^{n} &\le&  v_{xx}^{n-1} + {H(x_0, v_x^*,  q_*) - H(x_0, v_x^*,   v_{xx}^{n-1})\over  H_q(x_0, v_x^*,  v_{xx}^{n-1})} + C\e_n.
\eeaa
Denote $\dis\tilde \e_n := q_* - \inf_{m\ge n}  v_{xx}^{m-1}(x_0)\ge 0$. Then $\dis\lim_{n\to\infty}\tilde \e_n = 0$ and $q_* \le  v_{xx}^{n-1}(x_0) + \tilde \e_n$, thus
\beaa
 v_{xx}^{n} \le  v_{xx}^{n-1} - {1\over C_1} (  v_{xx}^{n-1} - q_*) + C(\e_n+\tilde \e_n).
\eeaa
This implies that
\beaa
 v_{xx}^{n} - q_* \le (1-{1\over C_1}) ( v_{xx}^{n-1} - q_*) + C(\e_n+\tilde \e_n).
\eeaa
Then by standard arguments we have $\limsup_{n\to\infty} ( v_{xx}^{n}(x_0) - q_*) \le 0$. This, together with the definition of $q_*$, implies the limit $\lim_{n\to\infty}  v_{xx}^{n}(x_0) = q_*$ exists.  Noting again that by \reff{scalar-paxvnest} $v^n_{xx}$ is uniformly bounded, then it follows from the closeness of the differentiation operator\footnote{Indeed, denoting $q^*(x) := \lim_{n\to\infty}  v_{xx}^{n}(x)$, then for any $x_1< x_2$ we have $v^n_x(x_2)-v^n_x(x_1) = \int_{x_1}^{x_2} v^n_{xx}(x)dx$. Sending $n\to \infty$, by the convergence of $v^n_x$ and the bounded convergence theorem, we have $v^*_x(x_2)-v^*_x(x_1) = \int_{x_1}^{x_2} q^*(x)dx$. This clearly implies that $v^*_{xx}(x) = q^*(x)$.}  
that $\lim_{n\to\infty}  v_{xx}^n(x) =  v_{xx}^*(x)$. Now send $n\to\infty$ in \reff{scalarPDEn}, we see that $v^*$ satisfies \reff{scalarHJBv}. Note further that $q\mapsto H(x,z,q)$ has an inverse function, then from \reff{scalarHJBv} we conclude that $ v_{xx}^*\in C^1_b(\dbR;\dbR)$. Finally, it follows from the uniqueness of classical solutions to \reff{scalarHJBv} that $v^*=v$. In particular, this implies that $ v\in C^3_b(\dbR;\dbR)$  and hence $\pi^*\in \cA_{Lip}$.
\qed

\subsection{Rate of convergence in a further special case}
In this subsection we obtain the rate of convergence under the following extra assumption.\footnote{\label{comparison}This condition can be viewed as a boundary condition (at $|x| = \infty$). This has quite different nature than the smallness condition in \cite{TWZ}, which roughly means $\si$ is not sensitive to the control $a$.  }

\begin{assum}
\label{assum-limit}
For $\f = b, \si, r$, there exist $\ul \f, \ol \f\in C^0_b(\dbR; \dbR)$ such that
\beaa
\lim_{x\to -\infty} \sup_{a\in A} |\f(x,a)- \ul \f(x)|=\lim_{x\to \infty} \sup_{a\in A} |\f(x,a)-\ol \f(x)|=0.
\eeaa
\end{assum}

\begin{thm}
\label{thm-small}
Let Assumptions \ref{assum-infinite2} and \ref{assum-limit} hold true. Then
\bea
\label{small-rate0}
\lim_{n\to\infty}\|\D v^n\|_2 + \|\D\pi^n\|_0 =0.
\eea
Moreover, there exist $\rho_0>0$, $C>0$,  and $0<\eta<1$, depending only on $d$ and $C_0$, such that, whenever $\rho\ge \rho_0$, \reff{small-rate0} has the following super-exponential rate: 
\bea
\label{small-rate}
\|\D v^n\|_2 + \|\D\pi^n\|_0 \le C \eta^{2^n}.
\eea
\end{thm}
\proof We proceed in three steps. Denote 
\bea
\label{en12}
\q \e^n_0 := \|\D v^n\|_0,\q \e^n_1:= \| \D v_x^n\|_0,\q \e^n_2:= \|  \D v_{xx}^n\|_0,\q \e^n_{1,2}:= \e^n_1 + \e^n_2.
\eea

{\bf Step 1.} Denote $\ol v(x) := {1\over \rho}[\ol r(x) + \l (\ln|A|)^+]$ and, for $R>0$,
\beaa
\d_R := \!\!\! \sup_{x\ge R, a\in A} \big[ |b(x,a)-\ol b(x)| + |\si(x,a)-\ol \si(x)|+|r(x,a)-\ol r(x)|\big]\to 0,
\eeaa
as $R\to\infty$. For any $x\ge 2R$, by \reff{scalarV} we have
\beaa
&&\dis \big|v(x) - \ol v(x)\big| \le \sup_{\pi\in \cA}\int_0^\infty e^{-\rho t} \dbE\big[|\tilde r(X^\pi_t, \pi(s, X^\pi_s)) - \ol r(x)|\big]dt\\
&&\dis \le {\d_R\over \rho} +  C\sup_{\pi\in \cA}\int_0^\infty e^{-\rho t} \dbP(X^\pi_t\le R)dt \le {\d_R\over \rho} +  C\sup_{\pi\in \cA}\int_0^\infty e^{-\rho t} \dbP(|X^\pi_t - x|\ge R)dt\\
&&\dis \le {\d_R\over \rho} +  {C\over R^2}\sup_{\pi\in \cA}\int_0^\infty e^{-\rho t} \dbE\big[|X^\pi_t - x|^2\big]dt \le {\d_R\over \rho} +  {C\over R^2}\int_0^\infty e^{-\rho t + C_1t} dt\\
&&\dis \le  {\d_R\over \rho} +  {C\over R^2(\rho-C_1)},
\eeaa
for some $C_1>0$, and here we assume $\rho\ge \rho_0>C_1$. Send $R\to \infty$,  clearly we have
\bea
\label{vinfty}
\lim_{x\to\infty} |v(x) - \ol v(x)| =0.
\eea

Next, by \reff{scalarHGamma}, \reff{scalarH}, and \reff{HpaqH2} one can easily show that 
\bea
\label{Hinfty}
\qq \sup_{x\ge R}\Big[ |H_z(\cd) - \ol b(x)| + |H_q(\cd)-{1\over 2}\ol \si^2(x)| + \big|h(\cd) + \ol v(x)\big|\Big](x,z,q)\le C_{z,q} \d_R,
\eea
 where $C_{z,q}$ depends on the bound of $z, q$. By \reff{scalarPDEn} we have
\beaa
&\dis v^n(x) = \dbE\Big[\int_0^\infty e^{-\rho t} h^{n-1}(X^n_t)dt\Big],\\
&\dis\mbox{where}\q X^n_t = x + \int_0^t H^{n-1}_z(X^n_s) ds + \int_0^t \sqrt{2H^{n-1}_q(X^n_s)}dW_s.
\eeaa
From the Step 1 in the proof of Theorem \ref{thm-main2}  we can easily see that, for $\rho$ large,  $v_x^{n-1}$ and $v_{xx}^{n-1}$ are uniformly bounded, uniformly in $n$ and $\rho$. Then, similarly to \reff{vinfty} we can show that $ \lim_{x\to\infty} \sup_n |v^n(x) - \ol v(x)| =0$, which in turn shows that   $ \lim_{x\to\infty} \sup_n |\D v^n(x)| =0$. Similarly we have $ \lim_{x\to -\infty} \sup_n |\D v^n(x)| =0$. These, together with Theorem \ref{thm-main2}, lead  to that 
\bea
\label{scalar-en0}
\lim_{n\to \infty} \e^n_0=0.
\eea

{\bf Step 2.}  Let $\rho_1>\rho$ be a large number. Recall \reff{scalar-rho1vnx}, similarly we have
 \bea
 \label{phiv2}
  v_x(x) := \dbE\Big[\int_0^\infty e^{-\rho_1 t} \big[f + (\rho_1-\rho) v\big] (X^x_t){W_t\over t}dt\Big],\nonumber\\
  \mbox{where}\q   f := H^v - {1\over 2} v_{xx},\q \phi^v (x) := \phi(x, v_x, v_{xx}).
 \eea
Here we abuse the notation $\phi^v$ with \reff{unrep} and \reff{vxnrep}. Then, 
\beaa
|\D v^n_x(x)| \le  \int_0^\infty e^{-\rho_1 t} \big[C+ (\rho_1-\rho) \e^n_0\big] {\dbE[|W_t|\over t}dt \le {C\over \sqrt{\rho_1}} +{C(\rho_1-\rho)\over \sqrt{\rho_1}} \e^n_0. 
\eeaa
By the arbitrariness of $x$, we have 
\beaa
\limsup_{n\to\infty} \e^n_1 \le {C\over \sqrt{\rho_1}} +{C(\rho_1-\rho)\over \sqrt{\rho_1}} \lim_{n\to \infty}\e^n_0={C\over \sqrt{\rho_1}}.
\eeaa
Since $\rho_1$ is arbitrary, we obtain 
\bea
\label{scalar-en1}
\lim_{n\to \infty} \e^n_1=0.
\eea

Moreover, recall \reff{scalarPDEn1} and similarly we have
\bea
\label{Hv}
\qq\q v_{xx} = {1\over   H^v_q }\big[ \rho v   - H^v_z  v_x-h^v \big]. 
 \eea
By \reff{Hinfty} and \reff{scalar-en0}, \reff{scalar-en1}, for $x\ge R$ we have
\beaa
|\D v_{xx}^n(x)| \le \Big|{2\over \ol \si^2 }\big[ \rho v^{n-1}   - \ol b v^{n-1}_x +\ol v \big] - {2\over \ol \si^2 }\big[ \rho v   - \ol b v_x +\ol v \big]\Big|(x) + C\d_R \le C\d_R.
 \eeaa
 That is, $\lim_{x\to\infty}\sup_n |\D v_{xx}^n(x)| =0$. Similarly $\lim_{x\to -\infty}\sup_n |\D v_{xx}^n(x)| =0$. Thus, it follows from Theorem \ref{thm-main2} that 
\bea
\label{scalar-en2}
\lim_{n\to \infty} \e^n_2=0.
\eea
Combining \reff{scalar-en0}, \reff{scalar-en1}, \reff{scalar-en2}, we obtain \reff{small-rate0}.

{\bf Step 3.} We now derive the rate of convergence when $\rho$ is large. Denote 
\beaa
X^x_t = x + \int_0^t b_0(X^x_s) ds+ \int_0^t \si_0(X^x_s) dW_s,~\mbox{where}~b_0 := H^v_z,~\si_0 := \sqrt{2H^v_q}.
\eeaa
Note that,  by \reff{scalarPDEn},
\bea
\label{Fn}
&\dis \rho \D v^n = {1\over 2}\si_0^2 \D v^n_{xx} + b_0 \D v^n_x+F^n,\q\mbox{where}\\
&\dis \dis F^n :=   H^{n-1}_q[ v_{xx}^{n}-  v_{xx}^{n-1}] + H^{n-1}_z[ v_x^{n} -  v_x^{n-1}]+ H^{n-1}- H^v - H_q^v \D v^n_{xx}  -H^v_z\D v^n_x.\nonumber
\eea
Since $H$ is jointly convex in $(z, q)$, we have
\beaa
0\le H^v  - H^{n-1} + H^{n-1}_q\D  v_{xx}^{n-1} + H^{n-1}_z\D  v_x^{n-1}  \le C\big[|\D  v_{xx}^{n-1}|^2 + |\D  v_x^{n-1}|^2\big].
\eeaa
Then, recalling the $\e^n_{12}$ in \reff{en12}, 
\bea
\label{Fnest}
|F^n|&\le&  \big|H^{n-1}_q- H^v_q\big||\D  v_{xx}^{n}| + \big|H^{n-1}_z- H^v_z\big||\D  v_x^{n}|+ C\big[|\D  v_{xx}^{n-1}|^2 + |\D  v_x^{n-1}|^2\big]\nonumber\\
\!\!\!&\le&\!\!\! C\Big[ [|\D v_x^{n-1}| + |\D v_{xx}^{n-1}|][ |\D v_{xx}^n| + |\D v_x^{n}|] +|\D  v_{xx}^{n-1}|^2 + |\D  v_x^{n-1}|^2\big]\nonumber\\
\!\!\!&\le&\!\!\! C\e^{n-1}_{1,2} \big[\e^n_{1,2} + \e^{n-1}_{1,2} \big].
\eea
Recall from Theorem \ref{thm-main2} that $v\in C^3_b(\dbR)$ and hence $ b_0, \si_0\in C^1_b(\dbR)$.  Then, for the $N^x_t$ corresponding to $b_0, \si_0$, we have
\beaa
&&\dis |\D v^n(x)| \le\hE\Big[\int_0^\infty e^{-\rho t} |F^n(X^x_t)| dt\Big] \le {C\over \rho}\e^{n-1}_{1,2} \big[\e^n_{1,2} + \e^{n-1}_{1,2} \big];\\ 
&&\dis |\D v_x^n(x)| \le \hE\Big[\int_0^\infty \!\!\!\!\! e^{-\rho t} |F^n(X^x_t) N^x_t |dt\Big] \le {C\over \sqrt{\rho-C_1}}\e^{n-1}_{1,2} \big[\e^n_{1,2} + \e^{n-1}_{1,2} \big];
\eeaa
for some appropriate $C_1$ and for $\rho\ge \rho_0> C_1+1$. Then we can easily get
\bea
\label{scalar-en01}
\rho \e^n_0 + \sqrt{\rho} \e^n_1 \le C\e^{n-1}_{1,2} \big[\e^n_{1,2} + \e^{n-1}_{1,2} \big].
\eea

Moreover, by \reff{Fn} and \reff{Fnest} we have
\beaa
\e^n_2 \le C\rho \e^n_0 + C\e^{n-1}_{1,2} \big[\e^n_{1,2} + \e^{n-1}_{1,2} \big]\le C\e^{n-1}_{1,2} \big[\e^n_{1,2} + \e^{n-1}_{1,2} \big].
\eeaa
Combined with \reff{scalar-en01}, this leads to 
$\e^n_{1,2}   \le C_2\e^{n-1}_{1,2} \big[\e^n_{1,2} + \e^{n-1}_{1,2} \big].$
By \reff{scalar-en1} and \reff{scalar-en2}, there exists $n_0$ such that $ \e^n_{1,2} \le {1\over 3C_2}$ for all $n\ge n_0$. 
Then 
\beaa
\left.\ba{c}
\e^n_{1,2} \le {3\over 2}C_2(\e^{n-1}_{1,2})^2,\q\mbox{and thus}\q  \e^n_{1,2} \le \big({3\over 2}C_2 \e^{n_0}_{1,2}\big)^{2^{n-n_0}} \le {1\over 2^{2^{n-n_0}}},\q n>n_0.
\ea\right.
\eeaa
This, together with \reff{scalar-en01},  implies
\bea
\label{scalar-ratev12}
\|\D v^n\|_2\le C \eta^{2n}, \q \text{for}~ \eta:= 2^{-2^{-n_0}} <1.
\eea
Here $C$ is chosen so that the above estimate holds true for $n\le n_0$ as well.
 
Finally, similarly to \reff{pirate} we can easily obtain the desired estimate for $\D \pi^n$.
\qed

\begin{rem}
\label{rem-smallrate}
Assumption \ref{assum-limit} is used to prove \reff{small-rate0}, but \reff{small-rate} relies only on \reff{small-rate0}, as we saw in Step 3 of the above proof. In other words, any possible alternative sufficient conditions for \reff{small-rate0} will also imply \reff{small-rate} when $\rho$ is large.
\end{rem}


\begin{thebibliography}{1} 
	
 \bibitem{BGMX}
Bai, L., Gamage, T. , Ma, J., and Xie, P., (2023) {\it Reinforcement Learning for Optimal Dividend Problem under Diffusion Model}, Preprint, arXiv:2309.10242. 

\bibitem{Bellman1}
Bellman, R. (1955) {\it Functional equations in the theory of dynamic programming}. {\sl V. Positivity and quasi-linearity, Proc. Nat. Acad. Sci. USA}, 41, 743-746.

\bibitem{Bellman2}
Bellman, R. (1957) {\sl Dynamic Programming}, Princeton University Press, Princeton, NJ.

\bibitem{Bismut}
Bismut, J. M., (1984) {\sl Large Deviation and Malliavin Calculus}, Progress in Mathematics {\bf 45}. Birkh\"a user.

\bibitem{BMZ}
Bokanowski, O., Maroso, S. and Zidani, H., (2009)  \emph{Some convergence results for Howard's
algorithm}. {\sl SIAM Journal on Numerical Analysis}, 47(4), 3001-3026.

\bibitem{Dong}
Dong, Y. (2024) {\it Randomized optimal stopping problem in continuous time and reinforcement learning algorithm}. {\sl SIAM J. Control Optim.}, 62 (3), 1590-1614.

\bibitem{EL}
Elworthy, K. D. and Li, X.M. (1994),  {\it Formulae for the derivatives of heat semigroups}. J. Funct. Anal. {\bf 125},  252--286. 

\bibitem{GXZ}
Guo, X.,  Xu, R., and  Zariphopoulou, T. (2022) {\it Entropy regularization for mean field games with learning}. {\sl Mathematics
of Operations research}, {\bf 47}(4), 3239--3260.

\bibitem{Howard}
Howard, R.A. (1960) {\sl Dynamic programming and Markov processes}, Cambridge MA, New York, The Technology Press of the MIT. J. Wiley.

\bibitem{Huang2023}
Huang, Y., Wang Z., and Zhou, Z. (2025),\emph{ Convergence of Policy Iteration for Entropy-Regularized Stochastic Control Problems}, {\sl SIAM J. Control Optim.} 63 (2), 752-777.



\bibitem{IRZ}
Ito, K., Reisinger, C. and Zhang, Y.,(2021),  \emph{A neural network-based policy iteration algorithm with global H2-superlinear convergence for stochastic games on domains. Foundations
of Computational Mathematics}, 21(2),331-374.

\bibitem{PIA1}
Jacka, S. and  Mijatovi\'{c}, A., (2017)  {\it On the policy improvement
  algorithm in continuous time}, {\sl Stochastics} \textbf{89}(1), 348--359. 

\bibitem{PIA2}
Kerimkulov, B., \v{S}i\v{s}ka, D., and Szpruch, L.,   (2020) \emph{Exponential
  convergence and stability of {H}oward's policy improvement algorithm for
  controlled diffusions}, {\sl SIAM J. Control Optim.} \textbf{58}(3),
  1314--1340. 
  
  
  \bibitem{PIA3}
Kerimkulov, B.,  \v{S}i\v{s}ka, D., and Szpruch, L.,  (2021) \emph{A modified {MSA} for
  stochastic control problems}, {\sl Appl. Math. Optim.} \textbf{84}(3),
  3417--3436. 

\bibitem{Krylov}
 Krylov, N. V. (1996) {\sl Lectures on elliptic and parabolic equations in H\"{o}lder spaces}, 
American Mathematical Society, Providence, RI.

\bibitem{MWZZ} 
Ma, J., Wang, G.,  Zhang, J., and Zhou, X.Y.,  {\it Reinforcement Learning Algorithms for Entropy-Regularized HJB Equations with Model Uncertainty}, working paper. 
  
\bibitem{MaZhang1} Ma, J. and Zhang, J. (2002),  {\it Representation Theorems for Backward Stochastic Differential Equations}. Annals of Applied Probability,  \textbf{12}(4), 1390--1418.   


\bibitem{Nualart}
Nualart,  D. (2006). {\sl Malliavin calculus and related topics}. In Stochastic Processes and Related Topics, 2nd ed., Springer, Berlin.

\bibitem{PIA4}
Puterman, M.~L.,  (1981)  \emph{On the convergence of policy iteration for controlled
  diffusions}, J. Optim. Theory Appl. \textbf{33}(1), 137--144.

\bibitem{PB}
 Puterman, M.~L., and  Brumelle, S.~L. (1979) {\it On the convergence of policy iteration in stationary dynamic programming}. {\sl Mathematics of Operations Research}, 4, 60-69.

\bibitem{RSZ}
Reisinger, C., Stockinger, W., and  Zhang, Y. (2023) {\it Linear convergence of a policy gradient
method for some finite horizon continuous time control problems}. {\sl SIAM Journal on
Control and Optimization}, 61, 3526-3558.

\bibitem{RZ}
Reisinger,  C. and Zhang, Y. (2021). {\it Regularity and stability of feedback relaxed controls}. {\sl SIAM J.
Control Optim.} {\bf 59}, 3118--3151.

\bibitem{SR}
Santos, M. S.  and Rust, J.,(2004), \emph{Convergence properties of policy iteration}. SIAM Journal on Control and Optimization, 42(6), 2094-2115.

\bibitem{SSZ}
Sethi, D., Siska, D., and Zhang, Y. (2024) {\it Entropy annealing for policy mirror descent in
continuous time and space}, Preprint, arXiv:2405.20250.

\bibitem{TZZ}
Tang, W.,  Zhang, Y. P., and Zhou, X. Y., (2022) {\it Exploratory HJB equations and their convergence}. {\sl SIAM Journal on
Control and Optimization}, {\bf 60}(6), 3191--3216.

\bibitem{TWZ}
Tran, H.~V., Wang, Z., and Zhang, Y.~.P., (2024) {\it Policy Iteration for Exploratory Hamilton-Jacobi-Bellman Equations}, 
Preprint, 	arXiv:2406.00612. 

\bibitem{WZZ}
Wang, H., Zariphopoulou, T. and Zhou, X.Y., (2020), \emph{Reinforcement
  learning in continuous time and space: a stochastic control approach}, {\sl J.
  Mach. Learn. Res}. \textbf{21}(198), 1--34.

\bibitem{WZ}
Wang, H. and  Zhou, X.Y.,  (2020), \emph{Continuous-time mean-variance portfolio
  selection: a reinforcement learning framework}, {\sl Math. Finance} \textbf{30}(4), 1273--1308. 

\bibitem{Zhangthesis}
Zhang, J. (2001), \emph{Some fine properties of backward stochastic differential equations}. Ph.D. dissertation, Purdue University.




\end{thebibliography}
\end{document}